# Théorème de Chebotarev et Congruences de suites récurrentes linéaires, liens avec les algorithmes de factorisations sur $\mathbb{F}_p$

par G.Duval

LMI, INSA Rouen France
guillaume.duval@insa-rouen.fr

**Résumé**

The classical congruences satisfied by the Fibonacci and Lucas sequences are reflected with the decompostion of primes in the ring generated by the gold number. This generalizes to establish a correspondence that we hope will be new between Chebotarev's theorem and the congruences satisfied by linear sequences. This link is done into the context of number field extensions. In particular we characterize primes ideals totally decomposed by simple congruences on the terms of linear recurrent sequences. Our results are illustrated by numerous examples, including Padovan sequencesof group $\mathfrak{S}_3$ and associated with the Cartier Trink polynomial of group $\mathbb{PSL}_2(\mathbb{F}_7)$. Furthermore, we establish a link between the factorisation algorithms of Berlekamp and Cantor-Zassenhaus and the results of this work.

**Résumé**

Les congruences classiques sur les suites de Fibonacci et Lucas se reflètent avec les décompositions des nombres premiers sur l'anneau engendré par le nombre d'or. Cela se généralise pour établir une correspondance que nous espérons nouvelle entre le théorème de Chebotarev et les congruences de suites linéaires. Ceci est fait dans le cas général des extensions de corps de nombres. En particulier nous caractérisons les idéaux premiers totalement décomposés par des congruences simples sur les termes de suites récurrentes linéaires. Nos résultats sont illustrés sur de nombreux exemples, dont les suites de Padovan de groupe $\mathfrak{S}_3$ et celles associées au polynôme de Cartier-Trinks de groupe $\mathbb{PSL}_2(\mathbb{F}_7)$. Par ailleurs, nous établissons un lien entre les algorithmes de factorisation de Berlekamp et de Cantor-Zassenhaus et les résultats de ce travail.

# 1 Introduction

Le point de départ de cette étude fut le constat du résultat suivant :

**Théorème 1.** *Soit $(U_n)_{n \in \mathbb{N}}$ une suite d'entiers relatifs vérifiant la récurrence de Fibonacci : $U_{n+2} = U_{n+1} + U_n$ alors, pour tout nombre premier $p$ les congruences $U_p \mod p$ sont très simples et données par la règle suivante :*

- *Si $p \equiv 1$ ou $4 \mod [5]$ alors $U_p \equiv U_1 \mod [p]$,*
- *Si $p \equiv 2$ ou $3 \mod [5]$ alors $U_p \equiv U_0 - U_1 \mod [p]$,*
- *Si $p = 5$ alors $U_5 = 3 U_0 + 5 U_1 \equiv 3 U_0 \mod [5]$.*

Ce théorème généralise les « congruences de Legendre » : $F_p \equiv 1 \mod [p]$ si $p \equiv 1$ ou $4 \mod [5]$ et $F_p \equiv -1 \mod [p]$ si $p \equiv 2$ ou $3 \mod [5]$ où $(F_n)_{n \in \mathbb{N}}$ est la suite de Fibonacci classique de conditions initiales $F_0 = 0, F_1 = 1$. Il englobe aussi les « congruences de Lucas » : $L_p \equiv 1 \mod [p]$, si $(L_n)_{n \in \mathbb{N}}$ vérifie la récurrence de Fibonacci avec les conditions initiales $L_0 = 2, L_1 = 1$.

Ce qui est plus intéressant est que ces congruences reflètent la manière dont les nombres premiers se factorisent dans le corps de nombres associé à la récurrence de Fibonacci. En effet, l'équation caractéristique de la récurrence de Fibonacci est $C(X) = X^2 - X - 1 = 0$. Si $E$ est le corps de décomposition de $C(X)$ sur $\mathbb{Q}$, alors pour tout nombre premier $p$ on a :

- $p$ est totalement décomposé dans $E$ si et seulement si $p \equiv 1$ ou $4 \mod [5]$, ou encore si $\left(\frac{5}{p}\right) = 1$.
- $p$ reste premier dans $E$ si et seulement si $p \equiv 2$ ou $3 \mod [5]$, ou encore si $\left(\frac{5}{p}\right) = -1$.
- le cas $p = 5$, correspondant au cas où $p$ divise le discriminant de $E$, c'est à dire quand $p$ est ramifié.

Nous allons voir que ce résultat illustre une correspondance entre d'un coté des congruence $U_p$ mod $[p]$ si $(U_n)_{n\in\mathbb{N}}$ est une suite d'entiers vérifiant une récurrence linéaire de polynôme caractéristique unitaire $C(X) = 0$ avec $C(X) \in \mathbb{Z}[X]$ et d'autre part la décomposition des idéaux premiers dans le corps de décomposition de $C(X)$.

Comme la théorie de Galois des extensions de corps de nombres et de leurs anneaux de Dedekind est une théorie relative, il préférable de se placer directement dans ce contexte général, même si dans les applications le corps de base sera majoritairement celui des rationnels. Nous adopterons donc le contexte général et les notations suivantes :

**Notation 2.** *Soient $E/K$ une extension galoisienne de corps de nombres, $A := \mathfrak{o}_K$ et $B := \mathfrak{o}_E$ les anneaux de Dedekind associés. Soit $C(X) \in A[X]$ un polynôme unitaire de degré $d \geqslant 2$ de racines simples $\{x_1, ..., x_d\}$ tel que $E$ soit le corps de décomposition de $C(X)$ sur $K$. Soient $\Delta_C$ et $\Delta_{B/A}$ les discriminants associés. Pour tout $\sigma$ dans le groupe de Galois $G(E/K)$ on notera $\mathbb{P}(E/K, \sigma) = \mathbb{P}(\sigma)$, l'ensemble des idéaux premiers $\mathcal{P}$ de $A$, premiers à $\Delta_C$[1], tels qu'il existe $Q \in \mathrm{Spec}(B)$ divisant $\mathcal{P}$ avec $\sigma = \mathrm{Frob}_Q$. Cet automorphisme de Frobenuis étant caractérisé par :*

$$\sigma(x) \equiv x^p \quad [Q], \tag{1}$$

*pour tout $x \in B$ si $p := N(\mathcal{P}) = N_\mathbb{Q}^K(\mathcal{P})$ est le cardinal du corps résiduel $\mathbb{F}_\mathcal{P} := A/\mathcal{P}$. Observons que le nombre $p$ est un nombre premier quand $K = \mathbb{Q}$, mais en général il coïncide avec la puissance d'un nombre premier. En particulier, $\mathbb{P}(E/K, \mathrm{Id}) := \mathrm{Split}(E/K)$ est l'ensemble des idéaux premiers de $K$ totalement décomposés dans $E$.*

**Définition 3.** *Soit $C(X) \in A[X]$ un polynôme unitaire. Dans toute la suite on dira qu'une suite $(U) = (U_n)_{n\in\mathbb{N}}$ est de classe $C(X)$ dans $A$ si*

- *Elle satisfait une récurrence linéaire d'équation caractéristique $C(X) = 0$.*
- *Pour tout $n \in \mathbb{N}$, $U_n \in A$.*

Afin de comprendre l'esprit de cet article à savoir cette correspondance entre les récurrences linéaires et la manière dont se décomposent les ideaux dans les corps de nombres, nous suggérons de commencer par consulter les énoncés suivants :

- Le Théorème 15 qui caractérise les idéaux totalement décomposés au moyen des restes $R_n(X)$, des divisions euclidiennes de $X^n$ par $C(X)$. Ce critère sera lui même illustré pour les suites de « Cartier Trinks » de groupe $\mathbb{P}\$\mathbb{L}_2(\mathbb{F}_7)$ en fin de ce texte.
- Le Théorème 18 qui généralise le Théorème 1 précédent pour toutes les extensions quadratiques. Comme conséquence de ce dernier résultat, le Corollaire 20 donne des congruences pour la fonction « tau » de Ramanujan.
- Pour les extensions Galoisiennes de plus petit groupe non abelien à savoir $\mathfrak{S}_3$ lesquelles correspondent à des récurrences linéaires d'ordre trois de type Padovan, le Théorème 25 et le Corollaire 28, illustrent encore une fois l'esprit de ce travail.

---

[1]. La condition $\mathcal{P}$ premier à $\Delta_C$ sert à garantir le fait que $\mathcal{P}$ ne se ramifie pas dans $B$. Par ailleurs, le discriminant $\Delta_B$ de $B$ est toujours un diviseur du discriminant $\Delta_C$ du polynôme $C(X)$. Enfin, le théorème Chebotarev garantit la non vacuité des ensembles $\mathbb{P}(\sigma)$, ces derniers étant infinis de densité $\mathrm{Card}(\mathcal{C}(\sigma))/\mathrm{Card}(G(E/K))$ où $\mathcal{C}(\sigma)$ est la classe de conjugaison de $\sigma$ dans $G(E/K)$.

- Enfin, dans la dernière section, nous établissons un lien entre les algorithmes de factorisation de Berlekamp et de Cantor-Zassenhaus et les résultats de ce travail.

Dans [5] Edouard Lucas utilise les notations symboliques suivantes que nous adopterons car elles sont commodes et suggestives : si $F(X) = \sum_k f_k X^k \in K[X]$ et si $(U_n)_{n \in \mathbb{N}}$ est une suite d'éléments de $A$ alors, nous noterons

$$F\{U_0\} := \sum_k f_k U_k.$$

Autrement dit, nous remplaçons les monômes $X^k$ par les termes $U_k$ correspondants de la suite. Pour obtenir un décalage nous écrirons $\sum_k f_k U_{k+n} := F\{U_n\} = (X^n \times F(X))\{U_0\}$.

Enfin, pour dire que $U = (U_n)_{n \in \mathbb{N}}$ est de classe $C$ dans $A$, nous écrirons $C\{U\} = 0$. Si $\deg(C) = d \geqslant 2$ et si $(U_n)_{n \in \mathbb{N}}$ est une suite de classe $C$ alors cette dernière dépend linéairement des conditions initiales $U_0, U_1, ..., U_{d-1}$. Par ailleurs, si $F(X)$ est de degré $\leqslant d-1$, la correspondance

$$U \mapsto F\{U_0\} = f_0 U_0 + \cdots + f_{d-1} U_{d-1},$$

sera une forme linéaire sur les suites de classes $C$.

Pour plus de détails sur ce calcul symbolique voir la Section 2.

Avec ces notations nous obtenons le résultat suivant :

**Théorème 4.** *Pour tout $\sigma \in G(E/K)$, il existe un polynôme $\Psi_\sigma(X)$ de degré inférieur à celui de $C(X)$ et à coefficients dans $\frac{1}{\Delta_C} B$ c'est à dire que $\Delta_C \Psi_\sigma(X) \in B[X]$ possédant les propriétés suivantes :*

1. *La correspondance $\sigma \mapsto \Psi_\sigma(X)$ de $G(E/K)$ vers le module $\frac{1}{\Delta_C} B[X]$ est injective*

2. *Si $\sigma \in Z(G(E/K))$, le centre du groupe de Galois de $E/K$, alors $\Delta_C \Psi_\sigma(X) \in A[X]$.*

3. *Si $\sigma \in Z(G(E/K))$, alors pour toute suite de classe $C$ dans $A$ et pour tout $\mathcal{P} \in \mathbb{P}(E/K, \sigma)$ on aura la congruence suivante dans $A$ :*

$$\Delta_C U_p = \Delta_C U_{N(\mathcal{P})} \equiv \Delta_C \Psi_\sigma\{U_0\} \quad [\mathcal{P}].$$

4. *En particulier pour $\sigma = \mathrm{Id}$, on aura $\Psi_{\mathrm{Id}}(X) = X$ et donc la congruence*

$$U_p \equiv U_1 \quad [\mathcal{P}], \quad \forall \mathcal{P} \in \mathrm{Split}(E/K).$$

5. *Pour tout $\sigma \in G(E/K)$, si $f \geqslant 1$ est l'ordre de $\sigma$ et si $\mathcal{P} \in \mathbb{P}(E/K, \sigma)$, alors la congruence précédente se généralise en*

$$U_{p^f} \equiv U_1 \quad [\mathcal{P}].$$

Une réciproque du Point 4, c'est à dire une caractérisation linéaire des idéaux complètement décomposés sera vue dans le Théorème 15 plus bas.

## 2 Calcul symbolique de Lucas avec les suites linéaires

Il est possible de démontrer le Théorème 4 sans faire référence à la présente section. Mais elle devient cependant indispensable pour la compréhension des illustrations de ce résultat. Par ailleurs nous introduirons les « suites de trace où de Lucas » qui joueront un rôle important par la suite.

### 2.1 Congruences polynomiales et égalités sur les suites

Soit $A$ un anneau et $C(X) \in A[X]$ un polynôme unitaire de degré $d \geqslant 2$. Soit $U = (U_n)_{n \in \mathbb{N}}$ une suite de classe $C$ dans $A$. Soit $R_n(X)$ le reste de la division Euclidienne de $X^n$ par $C(X)$. Ces polynômes sont tous de degrés $\leqslant d-1$ et à coefficients dans $A$ puisque $C(X)$ est unitaire. L'intérêt de cette suite est le suivant :

**Proposition 5.** *La suite des restes des divisions Euclidienne est de classe C dans $A[X]$. C'est à dire $C\{R\} = 0$ avec les notations symboliques de Lucas. Par ailleurs, pour toute suite $U = (U_n)_{n \in \mathbb{N}}$ de classe C dans A, chaque terme $U_n$ se calcule explicitement comme une forme linéaire en les conditions initiales $U_0, U_1, ..., U_{d-1}$ à coefficients dans A au moyen de la formule*

$$U_n = R_n\{U_0\}.$$

**Démonstration.** Écrivons $C(X) = c_0 + c_1 X + \cdots + c_{d-1} X^{d-1} + X^d$. Si $n$ et $k$ sont dans $\mathbb{N}$ on a la congruence $X^{n+k} \equiv R_{n+k}(X) \quad [C(X)]$ dans $A[X]$. On a alors les égalités et congruences suivantes

$$C\{R_n\} = c_0 R_n + c_1 R_{n+1} + \cdots + c_{d-1} R_{n+d-1} + R_{n+d}.$$

Ainsi $C\{R_n\} \equiv c_0 X^n + c_1 X^{n+1} + \cdots + c_{d-1} X^{n+d-1} + X^{n+d} = X^n C(X) \equiv 0 \quad [C(X)]$. Or le polynôme $C\{R_n\}$ est de degré $\leqslant d-1$ il est donc nul. Cela prouve donc que la suite des restes vérifie la récurrence $C\{R\} = 0$.

Si $U = (U_n)_{n \in \mathbb{N}}$ est une suite de classe C dans A, l'égalité $U_n = R_n\{U_0\}$ est vraie pour $n \leqslant d-1$ puisqu'alors $R_n(X) = X^n$. Supposons par récurrence que la relation $R_n\{U_0\} = U_n$ soit vraie jusqu'au rang $n \geqslant d-1$, on aura alors

$$\begin{aligned} U_{n+1} &= -c_{d-1} U_n - \cdots - c_1 U_{n+2-d} - c_0 U_{n+1-d} \\ &= -c_{d-1} R_n\{U_0\} - \cdots - c_1 R_{n+2-d}\{U_0\} - c_0 R_{n+1-d}\{U_0\} \\ &= R_{n+1}\{U_0\}, \end{aligned}$$

puisque $C\{R_{n+1-d}\} = 0$, ou encore que la suite des restes est de classe $C$. □

Cette propriété possède la généralisation suivante permettant d'obtenir des égalités de suites à partir de congruences polynomiales.

**Proposition 6.** *Soit $U = (U_n)_{n \in \mathbb{N}}$ est une suite de classe C dans A. Soit $F(X) \in A[X]$ si on a une division Euclidienne $F(X) = C(X) Q(X) + R(X)$ alors $F\{U_0\} = R\{U_0\}$. Plus généralement, si $F(X)$ et $G(X)$ dans $A[X]$ vérifient $F(X) \equiv G(X) \quad [C(X)]$ alors on aura $F\{U_0\} = G\{U_0\}$.*

**Démonstration.** Si nous écrivons $F(X) = \Sigma f_k X^k$ alors son reste s'exprimera par la formule : $R(X) = \Sigma f_k R_k(X)$. Ainsi nous aurons

$$F\{U_0\} = \Sigma f_k U_k = \Sigma f_k R_k\{U_0\} = R\{U_0\}.$$

Maintenant si $F(X) \equiv G(X) \quad [C(X)]$ alors $F$ et $G$ auront le même reste dans la division par $C$ et l'égalité $F\{U_0\} = G\{U_0\}$ résultera de la précédente. □

## 2.2 Les suites des traces ou de Lucas

Parmi toutes les suites de classes C dans l'anneau A, il en est une plus singulière que les autres à savoir la suite des traces ou de Lucas : Si $x_1, ..., x_d$ désignent encore les racines de $C(X)$, pour tout $n \in \mathbb{N}$, on posera :

$$L_n := x_1^n + \cdots + x_d^n.$$

Il y a une abondante littérature sur ce sujet par exemple [6]. Ce qui les rend singulières pour notre propos est la propriété suivante

**Proposition 7.** *Pour tout idéal premier $\mathcal{P}$ de A (ramifié ou non dans B) si $p = N(\mathcal{P})$ est le cardinal du corps résiduel $A/\mathcal{P}$. Alors pour tout $n \in \mathbb{N}$ on aura les congruences*

$$L_{p^n} \equiv L_1 \quad [\mathcal{P}].$$

Pour la preuve on se reportera encore à [6] où les auteurs soulignent le fait que dans le contexte où $A = \mathbb{Z}$, la propriété précédente n'est autre que l'analogue du théorème de Fermat matriciel. C'est à dire si $M \in \mathrm{GL}_n(\mathbb{Z})$ alors pour tout nombre premier on a la congruence

$$\mathrm{Tr}(M^p) \equiv \mathrm{Tr}(M) \quad [p].$$

Il en démontrent alors une très belle généralisation suivante pressentie par Arnold :

$$\text{Tr}(M^{p^{n+1}}) \equiv \text{Tr}(M^{p^n}) \quad [p^n].$$

## 2.3 Les relations universelles de Lucas

Soit $(U) = (U_n)_{n \in \mathbb{N}}$ une suite de classe $C(X)$ d'éléments de $A$. La relation $U_n = R_n\{U_0\}$, donne en particulier : si $p = N(\mathcal{P})$ est la norme d'un idéal premier quelconque de $A$ :

$$\begin{aligned} U_p &= R_p\{U_0\} \\ U_p &= \sum_{k=0}^{d-1} r_k(p)\, U_k, \end{aligned}$$

si le reste $R_p(X)$ s'exprime sous la forme :

$$R_p(X) = \sum_{k=0}^{d-1} r_k(p)\, X^k.$$

L'un des intérêts de la suite de Lucas et de la Proposition 7, est de nous fournir les congruences universelles suivantes :

**Proposition 8.** *Pour tout idéal premier $\mathcal{P}$ de $A$, (ramifié ou non dans $B$), si $p = N(\mathcal{P})$ alors on aura*

$$R_p\{L_0\} = \sum_{k=0}^{d-1} r_k(p)\, L_k \equiv L_1 \quad [\mathcal{P}].$$

**Démonstration.** $L_p = R_p\{L_0\} = \sum_{k=0}^{d-1} r_k(p)\, L_k \equiv L_1 \quad [\mathcal{P}]$ d'après la Proposition 7. □

Ici, nous avons utilisé l'expression « congruences universelles », car pour les suites de Lucas la congruence $U_p \bmod [\mathcal{P}]$ ne dépend pas de l'idéal premier $\mathcal{P}$, ce qui n'est pas le cas pour les autres suites de classe $C$.

## 3 Preuve du Théorème 4 et d'une réciproque partielle

Cette preuve se fera en plusieurs étapes.

### 3.1 Une première congruence liées aux Frobenuis

Soit $(U) = (U_n)_{n \in \mathbb{N}}$ de classe $C$ dans $A$. Puisque les racines $\{x_1, ..., x_d\}$ de $C(X)$ sont simples, la suite $(U)$ peut s'écrire comme combinaison linéaire des suites géométriques de raisons $x_i$ sous la forme

$$U_n = \lambda_1 x_1^n + \cdots + \lambda_d x_d^n, \tag{2}$$

où les coefficients $\lambda_i$ sont dans $E$ et ne dépendent pas de $n$. En écrivant les $d$ premières relations correspondantes pour $0 \leqslant n \leqslant d-1$, on obtient une relation matricielle liant les coefficients $\lambda_i$ aux conditions initiales de la suite :

$$\begin{pmatrix} U_0 \\ U_1 \\ \vdots \\ U_{d-1} \end{pmatrix} = V \begin{pmatrix} \lambda_1 \\ \vdots \\ \lambda_d \end{pmatrix} = \begin{pmatrix} 1 & \cdots & 1 \\ x_1 & \cdots & x_d \\ \vdots & & \vdots \\ x_1^{d-1} & \ldots & x_d^{d-1} \end{pmatrix} \begin{pmatrix} \lambda_1 \\ \vdots \\ \lambda_d \end{pmatrix}, \tag{3}$$

où $V := V(x_1, ..., x_d)$ est la matrice de Vandermonde de la famille $\{x_1, ..., x_d\}$. Cette matrice est de taille $d \times d$ à coefficients dans $B$ et son déterminant vaut

$$\det(V) := \prod_{d \geqslant i > j \geqslant 1} (x_i - x_j) := \sqrt{\Delta_C}. \tag{4}$$

Ainsi, en multipliant (3) par la transposée de la comatrice on obtient tout l'information importante suivante : chaque $\lambda_i$ appartient à $\frac{1}{\sqrt{\Delta_C}} B \subset \frac{1}{\Delta_C} B$ car $\sqrt{\Delta_C} \in B$.

**Lemme 9.** *Soit $\sigma \in G(E/K)$ et $Q \in \operatorname{Spec}(B)$ tel que $\sigma = \operatorname{Frob}(Q)$. Soit $\mathcal{P} := Q \cap A$ et $p = N(\mathcal{P})$. Soit $(U)$ une suite de classe $C$ se décomposant sous la forme donnée par (2). Alors la permutation des racines de $C(X)$ donnée par l'automorphisme de Frobenius $\sigma = \operatorname{Frob}(Q)$ donne la congruence :*

$$\Delta_C U_p \equiv \Delta_C [\lambda_1 \sigma(x_1) + \cdots + \lambda_d \sigma(x_d)] \quad [Q], \tag{5}$$

**Démonstration.** Si $n = N(\mathcal{P}) = p$ et si $\sigma = \operatorname{Frob}(Q)$, on aura

$$\begin{aligned} U_p &= \lambda_1 x_1^p + \cdots + \lambda_d x_d^p. \\ \Delta_C U_p &= (\Delta_C \lambda_1) x_1^p + \cdots + (\Delta_C \lambda_d) x_d^p \\ \Delta_C U_p &\equiv (\Delta_C \lambda_1) \sigma(x_1) + \cdots + (\Delta_C \lambda_d) \sigma(x_d) \quad [Q] \\ \Delta_C U_p &\equiv \Delta_C [\lambda_1 \sigma(x_1) + \cdots + \lambda_d \sigma(x_d)] \quad [Q] \end{aligned}$$

puisque chaque racine $x_i$ est dans $B$ de même que les nombres $\Delta_C \lambda_i$ et qu'enfin l'automorphisme de Frobenuis $\sigma = \operatorname{Frob}(Q)$ est caractérisé par les congruences (1) : $\sigma(x) \equiv x^p \quad [Q]$. □

**Remarque 10.** Ce lemme suggère les deux remarques suivantes :

1. Dans la congruence $\Delta_C U_p \equiv \Delta_C [\lambda_1 \sigma(x_1) + \cdots + \lambda_d \sigma(x_d)] \quad [Q]$, l'intérêt de l'expression $\Delta_C [\lambda_1 \sigma(x_1) + \cdots + \lambda_d \sigma(x_d)]$ est de ne dépendre que de la permutation induite sur les racines $x_i$ et plus du tout des idéaux premiers. Ces derniers sont en nombre infini alors que nous n'avons qu'un nombre fini de permutations des racines. Cela va être tout l'esprit de cet article que d'exhiber des congruences du type :

$$U_p \equiv y \quad [Q],$$

où les quantités $y$ appartiennent à un ensemble fini de combinaisons linéaires des conditions initiales $(U_0, ..., U_{d-1})$ et plus des idéaux premiers $Q$ et $Q \cap A$.

2. Pour chaque $\sigma \in G(E/K)$, l'expression $\lambda_1 \sigma(x_1) + \cdots + \lambda_d \sigma(x_d)$ devient linéaire en les conditions initiales $U_0, ..., U_{d-1}$ après inversion de (3). C'est le calcul que nous allons faire maintenant.

## 3.2 Les polynômes permutants des racines

D'après [4] p. 287 on a le résultat qui peut être modifié sous la forme équivalente suivante :

**Lemme 11.** *Soit $C(X) = \prod_{i=1}^{d}(X - x_i)$ un polynôme unitaire à racines simples dans un corps $E$. Si nous posons :*

$$\frac{C(X)}{(X - x_i)} := \beta_0(x_i) + \beta_1(x_i) X + \cdots + \beta_{d-1}(x_i) X^{d-1}$$

*Alors pour tout entier $0 \leqslant r \leqslant d-1$. On a l'égalité polynomiale :*

$$X^r = \sum_{i=1}^{d} \frac{C(X)}{(X - x_i)} \frac{x_i^r}{C'(x_i)}.$$

*L'inverse de la matrice de Vandermonde est donnée par*

$$(V(x_1, ..., x_d))^{-1} = \begin{pmatrix} \frac{\beta_0(x_1)}{C'(x_1)} & \frac{\beta_1(x_1)}{C'(x_1)} & \cdots & \frac{\beta_{d-1}(x_1)}{C'(x_1)} \\ \frac{\beta_0(x_2)}{C'(x_2)} & \frac{\beta_1(x_2)}{C'(x_2)} & \cdots & \frac{\beta_{d-1}(x_2)}{C'(x_2)} \\ \vdots & \vdots & & \vdots \\ \frac{\beta_0(x_d)}{C'(x_d)} & \frac{\beta_1(x_d)}{C'(x_d)} & \cdots & \frac{\beta_{d-1}(x_d)}{C'(x_d)} \end{pmatrix}.$$

Les calculs de ce résultat classique étant utiles pour la suite, nous en donnons la preuve.

**Démonstration.** Soit $r$ un entier naturel compris entre 0 et $d-1$. Le polynôme $\sum_{i=1}^{d}\frac{C(X)}{(X-x_i)}\frac{x_i^r}{C'(x_i)}$ est de degré $\leqslant d-1$ et vaut $x_i^r$ quand on l'évalue en $X = x_i$, pour tout $1 \leqslant i \leqslant d$. Il coïncide donc avec $X^r$. On a donc les égalités

$$\begin{aligned} X^r &= \sum_{i=1}^{d} \frac{C(X)}{(X-x_i)}\frac{x_i^r}{C'(x_i)}. \\ &= \sum_{i=1}^{d}\left(\sum_{k=0}^{d-1}\frac{\beta_k(x_i)}{C'(x_i)}x_i^r X^k\right) \\ X^r &= \sum_{k=0}^{d-1}\left(\sum_{i=1}^{d} x_i^r \times \frac{\beta_k(x_i)}{C'(x_i)}\right)X^k. \end{aligned}$$

Ainsi pour tous entiers $r$ et $k$ compris entre 0 et $d-1$, la somme $\sum_{i=1}^{d} x_i^r \times \frac{\beta_k(x_i)}{C'(x_i)}$ est égale au symbole de Kronecker $\delta_{r,k}$. Cela signifie que le produit de $V(x_1,...,x_d)$ et de la matrice du lemme coïncide avec la matrice identité de taille $d$. $\square$

**Définition 12.** *Pour toute permutation $\sigma \in G(E/K)$ posons le vecteur ligne :*

$$(\Psi_0(\sigma),...,\Psi_{d-1}(\sigma)) := (\sigma(x_1),...,\sigma(x_d)) \times V^{-1}$$

*et associons lui le polynôme :* $\Psi_\sigma(X) := \Psi_0(\sigma) + \cdots + \Psi_{d-1}(\sigma) X^{d-1}$

Avec ces notations nous obtenons le :

**Lemme 13.** *L'expression $\lambda_1 \sigma(x_1) + \cdots + \lambda_d \sigma(x_d)$ de (5) s'exprime en fonction des conditions initiales sous la forme :*

$$\lambda_1\sigma(x_1) + \cdots + \lambda_d\sigma(x_d) = \Psi_0(\sigma) U_0 + \cdots + \Psi_{d-1}(\sigma) U_{d-1} := \Psi_\sigma\{U_0\}.$$

*Pour tout $0 \leqslant k \leqslant d-1$, le coefficient $\Psi_k(\sigma)$ s'exprime par la formule*

$$\Psi_k(\sigma) = \sum_{i=0}^{d} \frac{\beta_k(x_i)}{C'(x_i)} \sigma(x_i).$$

*Plus synthétiquement*

$$\Psi_\sigma(X) = \sum_{i=1}^{d} \frac{C(X)}{(X-x_i)} \frac{\sigma(x_i)}{C'(x_i)}.$$

*Pour toute racine $x_i$ de $C(X)$, $\Psi_\sigma(x_i) = \sigma(x_i)$. Enfin, les polynômes $\Delta_C \Psi_\sigma(X)$ appartiennent à $B[X]$.*

La condition $\Psi_\sigma(x_i) = \sigma(x_i)$ pour toute racine $x_i$ de $C(X)$ et le fait que $\deg(\Psi_\sigma(X)) \leqslant d-1$ caractérise $\Psi_\sigma$ comme étant l'unique polynôme d'interpolation de Lagrange qui permute les racines de $C(X)$ selon la permutation $\sigma$. D'où le fait que dans le titre de cette section nous ayons nommé les $\Psi_\sigma$ : « **polynômes permutants des racines** ».

**Démonstration.** D'après (3) nous avons

$$\begin{aligned} \lambda_1\sigma(x_1) + \cdots + \lambda_d\sigma(x_d) &= (\sigma(x_1),...,\sigma(x_d)) \times (\lambda_1,...,\lambda_d)^T \\ &= (\sigma(x_1),...,\sigma(x_d)) \times V^{-1} \times (U_0,...,U_{d-1})^T \\ &= (\Psi_0(\sigma),...,\Psi_{d-1}(\sigma)) \times (U_0,...,U_{d-1})^T \\ \lambda_1\sigma(x_1) + \cdots + \lambda_d\sigma(x_d) &= \Psi_0(\sigma) U_0 + \cdots + \Psi_{d-1}(\sigma) U_{d-1} \\ &= \Psi_\sigma\{U_0\}, \end{aligned}$$

avec les notations symboliques de Lucas. L'expression $\Psi_k(\sigma) = \sum_{i=0}^{d} \frac{\beta_k(x_i)}{C'(x_i)} \sigma(x_i)$ s'obtient en faisant le produit matriciel $(\sigma(x_1), ..., \sigma(x_d)) \times V^{-1}$ avec la formule explicite de $V^{-1}$ donnée par le Lemme 11. Ainsi,

$$\Psi_k(\sigma) = \sum_{i=0}^{d} \sigma(x_i) \frac{\beta_k(x_i)}{C'(x_i)}$$

$$\sum_{k=0}^{d-1} \Psi_k(\sigma) X^k = \sum_{k=0}^{d-1} \left( \sum_{i=0}^{d} \frac{\beta_k(x_i)}{C'(x_i)} \sigma(x_i) \right) X^k$$

$$= \sum_{i=1}^{d} \left( \sum_{k=0}^{d-1} \frac{\beta_k(x_i)}{C'(x_i)} X^k \right) \sigma(x_i)$$

$$\Psi_\sigma(X) = \sum_{i=1}^{d} \frac{C(X)}{(X-x_i)} \frac{\sigma(x_i)}{C'(x_i)},$$

toujours d'après le Lemme 11 précédent.

Enfin si dans cette formule nous remplaçons $X$ par $x_i$, nous obtenons $\Psi_\sigma(x_i) = \sigma(x_i)$.

Les polynômes $\Psi_\sigma(X)$ sont à coefficients dans $E$ plus précisément d'après (4),

$$\sqrt{\Delta_C} \, \Psi_\sigma(X) \in B[X] \Rightarrow \Delta_C \, \Psi_\sigma(X) \in B[X].$$

□

Si on étend l'action de $G(E/K)$ à $E(X)$ en faisant agir les automorphismes sur les coefficients des fractions en laissant $X$ invariant nous obtenons alors le

**Lemme 14.** *L'application de $G(E/K)$ dans $E[X]$ associant $\sigma \mapsto \Psi_\sigma(X)$ est injective.*

*Pour tous $\theta$ et $\sigma$ dans $G(E/K)$, on a la relation de conjugaison :*

$$\theta(\Psi_\sigma(X)) = \Psi_{\theta\sigma\theta^{-1}}(X).$$

*Ainsi, $\Psi_\sigma(X) \in K[X]$ si et seulement si, $\sigma$ est dans le centre de $G$.*

**Démonstration.** Soient $\sigma$ et $\sigma'$ dans $G(E/K)$ si $\Psi_\sigma(X) = \Psi_{\sigma'}(X)$. Alors pour toutes les racines $x_i$ de $C(X)$, nous aurons $\sigma(x_i) = \Psi_\sigma(x_i) = \Psi_{\sigma'}(x_i) = \sigma'(x_i)$ et comme l'extension $E/K$ est engendrée par les racines de $C(X)$, nous aurons $\sigma = \sigma'$. Ceci prouve que la correspondance $\sigma \mapsto \Psi_\sigma(X)$ est injective.

Soit $\theta$ dans $G(E/K)$, distinct ou non de $\sigma$ nous aurons

$$\theta(\Psi_\sigma(X)) = \sum_{i=1}^{d} \frac{C(X)}{(X-\theta(x_i))} \frac{\theta(\sigma(x_i))}{C'(\theta(x_i))}$$

$$= \sum_{i=1}^{d} \frac{C(X)}{(X-\theta(x_i))} \frac{\theta\sigma\theta^{-1}(\theta(x_i))}{C'(\theta(x_i))}$$

$$= \sum_{l=1}^{d} \frac{C(X)}{(X-x_l)} \frac{\theta\sigma\theta^{-1}(x_l)}{C'(x_l)}$$

$$\theta(\Psi_\sigma(X)) = \Psi_{\theta\circ\sigma\circ\theta^{-1}}(X).$$

Pour passer de la ligne 2 à la ligne 3 des équations précédentes, nous avons juste numéroté les racines avec la permutation $x_l := \theta(x_i)$.

Nous avons donc les équivalences suivantes $\Psi_\sigma(X) \in K[X]$ si et seulement si pour tout $\theta \in G(E/K)$, $\theta(\Psi_\sigma(X)) = \Psi_\sigma(X) = \Psi_{\theta\circ\sigma\circ\theta^{-1}}(X)$. Et ceci se produit si et seulement si

$$\sigma = \theta \circ \sigma \circ \theta^{-1} \quad \forall \theta \in G(E/K),$$

car la correspondance $s \mapsto \Psi_s(X)$ est injective. Et cette dernière équivalence a lieu si et seulement si $\sigma$ est dans le centre de $G(E/K)$. □

### 3.3 Preuve du Théorème 4

Les Points 1 et 2 du théorème sont contenus dans les Lemmes 13 et 14.

#### 3.3.1 Preuve du Point 3.

Soit $\sigma$ dans $G(E/K)$, $\mathcal{P} \in \mathbb{P}(\sigma)$ et $Q \in \operatorname{Spec}(B)$ au dessus de $\mathcal{P}$ tel que $\operatorname{Frob}_Q = \sigma$. D'après la congruence (5) et le Lemme 13, nous avons

$$\Delta_C U_p \equiv \Delta_C \Psi_\sigma\{U_0\} \quad [Q]. \tag{6}$$

Les deux membres de cette congruence sont dans $B$. Maintenant, si $\sigma$ est dans le centre de $G(E/K)$, les deux membres de (6) sont dans $A$, d'après le Lemme 14. Puisque $Q \cap A = \mathcal{P}$, nous en déduisons que

$$\Delta_C U_p \equiv \Delta_C \Psi_\sigma\{U_0\} \quad [\mathcal{P}]$$

#### 3.3.2 Preuve du Point 4.

Enfin, quand $\sigma = \operatorname{Id}$, c'est à dire quand $\mathcal{P} \in \operatorname{Split}(E/K)$, le polynôme $\Psi_\sigma$ induit la permutation identité sur les racines $x_i$ on a donc $\Psi_{\operatorname{Id}}(X) = X$. C'est ce que prévoyaient déjà les Lemmes 13 et 11. Et ces derniers impliquent aussi que nous ayons :

$$\Psi_{\operatorname{Id}}(X) = \sum_{i=1}^d \frac{C(X)}{(X-x_i)} \frac{x_i}{C'(x_i)} = X \Rightarrow \Psi_{\operatorname{Id}}\{U_0\} = U_1,$$

Ceci donne comme précédemment $\Delta_C (U_p - U_1) \in \mathcal{P}$, mais puisque $\Delta_C \notin \mathcal{P}$ et que $\mathcal{P}$ est premier, nous obtenons finalement $U_p \equiv U_1 \quad [\mathcal{P}]$ quand $\mathcal{P} \in \operatorname{Split}(E/K)$.

#### 3.3.3 Preuve du Point 5.

Si $\sigma$ est d'ordre $f \geqslant 1$ dans le groupe de Galois, soit $\mathcal{P} \in \mathbb{P}(E/K, \sigma)$ et $Q$ divisant $\mathcal{P}$ dans $B$ tel que $\sigma = \operatorname{Frob}(Q)$. Puisque $\mathcal{P}$ n'est pas ramifié, $f = f_\mathcal{P}$ est le degré résiduel de l'extension $[\mathbb{F}_Q/\mathbb{F}_\mathcal{P}] = [B/Q : A/\mathcal{P}]$. Ainsi pour tout $x \in B$, nous aurons :

$$x = \sigma^f(x) \equiv x^{p^f}, \quad [Q].$$

Donc toujours d'après la relation (2), nous aurons :

$$\begin{aligned}
U_{p^f} &= \lambda_1 x_1^{p^f} + \cdots + \lambda_d x_d^{p^f}, \\
\Delta_C U_{p^f} &= (\Delta_C \lambda_1) x_1^{p^f} + \cdots + (\Delta_C \lambda_d) x_d^{p^f} \\
\Delta_C U_{p^f} &\equiv (\Delta_C \lambda_1) x_1 + \cdots + (\Delta_C \lambda_d) x_d, \quad [Q] \\
\Delta_C U_{p^f} &\equiv \Delta_C U_1 \quad [Q].
\end{aligned}$$

Cette dernière congruence signifie que $\Delta_C(U_{p^f} - U_1) \in Q \cap A = \mathcal{P}$. Comme $\Delta_C \notin \mathcal{P}$, on a bien que $U_{p^f} - U_1 \in \mathcal{P}$.

### 3.4 Caractérisation linéaire des idéaux totalement décomposés

Avec les notations habituelles, le Théorème 4, donne l'implication :

$$\mathcal{P} \in \operatorname{Split}(E/K) \Rightarrow U_p \equiv U_1 \quad [\mathcal{P}],$$

pour toute suite $(U)$ de classe $C(X)$ de $A$. Parmi toutes les suites de classe $C(X)$, la Proposition 5, signifie que la suite des restes $(R_n(X))_{n \in \mathbb{N}}$ est plus importante que les autres. Et la relation $U_p = R_p\{U_0\}$ permet de voir aisément en jouant avec les conditions initiales que la condition $U_p \equiv U_1 \quad [\mathcal{P}]$ pour toute suite de classe $C(X)$ d'éléments de $A$, équivaut à la condition plus synthétique :

$$R_p(X) \equiv X \quad [\mathcal{P}].$$

**Théorème 15.** *Soit $R_n(X)$, le reste de la division Euclidienne de $X^n$ par $C(X)$. Soit $\mathcal{P}$ un idéal premier de $A$ non ramifié dans $B$. Alors $\mathcal{P}$ est totalement décomposé dans $B$ si et seulement si $R_p(X) \equiv X \quad [\mathcal{P}]$ où $p = N(\mathcal{P})$.*

Nous aurons besoin du Lemme suivant qui doit être classique, ainsi que sa preuve...

**Lemme 16.** *Soit $C(X) \in K[X]$ à racines simples $x_1, ..., x_d$ dans son corps de décomposition $E$. Soit $\mathcal{P}$ un idéal premier de $A := \mathfrak{o}_K$ non ramifié dans $B := \mathfrak{o}_E$ c'est à dire que $\mathcal{P}$ est premier au discriminant $\Delta_C$. Si $Q$ divise $\mathcal{P}$ dans $B$ alors on a*
$$\mathbb{F}_Q = \mathbb{F}_\mathcal{P}[\bar{x}_1, ..., \bar{x}_d],$$
*si $\mathbb{F}_Q := B/Q$ et $\mathbb{F}_\mathcal{P} := A/\mathcal{P}$ sont les corps résiduels associés.*

**Démonstration. du Lemme** Soit $\mathbb{F} := \mathbb{F}_\mathcal{P}[\bar{x}_1, ..., \bar{x}_d]$, vu comme sous corps intermédiaire de la chaîne $\mathbb{F}_\mathcal{P} \subset \mathbb{F} \subset \mathbb{F}_Q$. L'extension $\mathbb{F}/\mathbb{F}_\mathcal{P}$ est normale donc galoisienne. Soient les deux morphismes de réductions : $\mathcal{R}_1: G_Q \to \mathrm{Gal}(\mathbb{F}_Q/\mathbb{F}_\mathcal{P})$ où $G_Q$ est le sous groupe de décomposition de $Q$ et $\mathcal{R}_2: \mathrm{Gal}(\mathbb{F}_Q/\mathbb{F}_\mathcal{P}) \to \mathrm{Gal}(\mathbb{F}/\mathbb{F}_\mathcal{P})$. Puisque $\mathcal{P}$ n'est pas ramifié, $\mathcal{R}_1$ est un isomorphisme. Par ailleurs, $\mathcal{R}_2$ est surjectif d'après le théorème de correspondance de Galois. Il suffit de voir que $\mathcal{R}_2 \circ \mathcal{R}_1$ est injectif pour déduire que $\mathcal{R}_2$ est lui aussi un isomorphisme et conclure que $\mathbb{F} = \mathbb{F}_Q$, toujours d'après le théorème de correspondance de Galois.

Soit $s \in \mathrm{Ker}(\mathcal{R}_2 \circ \mathcal{R}_1)$. On aura en particulier $s(\bar{x}_i) = \bar{x}_i$. Ce qui implique $s(x_i) \equiv x_i[Q]$, pour toute racine $x_i$. Si $s$ ne coïncidait pas avec l'identité de $E$, il existerait une racine $x_i$ telle que $s(x_i) \neq x_i$. Qui à renuméroter, nous aurions par exemple $s(x_1) = x_2$. De $\Delta_C = \Pi_{i<j}(x_i - x_j)^2$, le discriminant se factorise dans $B$ par $(x_1 - x_2)^2$. Mais $(x_1 - x_2)^2 = (x_1 - s(x_1))^2$ appartiendrait à $Q$ donc $\Delta_C$ aussi. Nous aurions donc que $\Delta_C$ appartiendrait à $Q \cap A = \mathcal{P}$ ce qui est contraire à l'hypothèse. Finalement, $s = \mathrm{Id}$ et $\mathcal{R}_2 \circ \mathcal{R}_1$ est injectif. $\square$

**Démonstration. du Théorème**. L'implication : si $\mathcal{P} \in \mathrm{Split}(E/K)$ alors $R_p(X) \equiv X[\mathcal{P}]$ est une conséquence du Théorème 4 et fut vue plus haut. Réciproquement, soit $\mathcal{P}$ non ramifié tel que $R_p(X) \equiv X[\mathcal{P}]$ si $p = N(\mathcal{P})$. Soit $Q$ divisant $\mathcal{P}$ dans $B$ et $\sigma := \mathrm{Frob}_Q$. Pour toute racine on aura $\sigma(x_i) \equiv x_i^p[Q]$. Mais $x_i^p = R_p(x_i)$ et par hypothèse $R_p(x_i) \equiv x_i[\mathcal{P}B]$. Donc nous aurons $x_i^p \equiv x_i[Q]$ puis $\sigma(x_i) \equiv x_i[Q]$ puisque $\mathcal{P}B \subset Q$. Puisque $\mathbb{F}_Q = \mathbb{F}_\mathcal{P}[\bar{x}_1, ..., \bar{x}_d]$, d'après le lemme précédent, $\sigma$ agit comme l'identité sur le corps résiduel $\mathbb{F}_Q$. Mais $\mathrm{Gal}(\mathbb{F}_Q/\mathbb{F}_\mathcal{P})$ est engendré par le Frobenuis, donc $\mathbb{F}_Q = \mathbb{F}_\mathcal{P}$ ce qui prouve que $\mathcal{P}$ est totalement décomposé. $\square$

# 4 Quelques illustrations du Théorème 4

## 4.1 Les cas abeliens et quadratiques

La plus immédiate est le

**Corollaire 17.** *Si l'extension de corps de nombres $E/K$ est abelienne alors pour tout $\sigma \in G(E/K)$, il existe un polynôme $\Psi_\sigma(X)$ appartenant au module $\frac{1}{\Delta_C}A[X]$ telle que pour toute suite de classe $C$ dans $A$ et pour tout $\mathcal{P} \in \mathbb{P}(E/K, \sigma)$ on ait la congruence suivante dans $A$ :*
$$\Delta_C U_p = \Delta_C U_{N(\mathcal{P})} \equiv \Delta_C \Psi_\sigma\{U_0\} \quad [\mathcal{P}].$$

Les choses sont toujours plus simples dans le cas abelien. Avant d'envisager des contextes plus complexes, le résultat suivant pour les suites linéaires d'ordre deux, englobe le cas des suites de Fibonacci évoqué dans le Théorème 1.

**Théorème 18.** *Soit $C(X) = X^2 - sX + \pi \in A[X]$ un polynôme de discriminant $\Delta_C = s^2 - 4\pi$ non carré dans $A$ de sorte que l'extension $E/K$ soit quadratique de groupe de Galois $G = <\tau> \simeq \mathbb{Z}/2\mathbb{Z}$. Alors pour toute suite $(U_n)_{n \in \mathbb{N}}$ de classe $C(X)$ dans $A$. Nous avons les congruences suivantes :*

- *Si $\mathcal{P} \in \mathrm{Split}(E/K)$ alors $U_p = U_{N(\mathcal{P})} \equiv U_1 \quad [\mathcal{P}].$*

- Si $\mathcal{P} \in \mathbb{P}(\tau)$ alors $U_p = U_{N(\mathcal{P})} \equiv s\, U_0 - U_1 \quad [\mathcal{P}]$.

Bien entendu quand $A = \mathbb{Z}$, les nombres premiers appartenant à $\text{Split}(E/\mathbb{Q})$ sont ceux pour lesquels $\left(\frac{\Delta_C}{p}\right) = 1$ et le Théorème de Réciprocité Quadratique permet de les identifier simplement.

**Démonstration.** La première congruence fut vue précédemment pour la seconde, soit $\tau$ la transposition échangeant les deux racines $x_1$ et $x_2$ de $C(X)$. Le polynôme transposant les deux racines est alors :

$$\begin{aligned}
\Psi_\tau(X) &= \frac{C(X)}{(X-x_1)} \frac{x_2}{C'(x_1)} + \frac{C(X)}{(X-x_2)} \frac{x_1}{C'(x_2)} \\
&= (X-x_2)\frac{x_2}{(x_1-x_2)} - (X-x_1)\frac{x_1}{(x_1-x_2)} \\
&= (x_1+x_2) - X \\
\Psi_\tau(X) &= s - X.
\end{aligned}$$

(Cette égalité pouvait être vue plus vite en observant que $\Psi_\tau(X)$ est le polynôme de degré un échangeant $x_1$ et $x_2$). Si $\mathcal{P} \in \mathbb{P}(\tau)$ alors $U_p = U_{N(\mathcal{P})} \equiv \Psi_\tau\{U_0\} = s\, U_0 - U_1 \quad [\mathcal{P}]$, avec les notations symboliques de Lucas. $\square$

### 4.2 Un résultat d'unicité pour les récurrences d'ordre deux

Il s'agit de la conséquence suivante du Théorème 18 :

**Corollaire 19.** *Soit $C(X) = X^2 - s\, X + \pi \in A[X]$ avec $s \neq 0$. Soient $(U)$ et $(V)$ deux suites de classe $C(X)$ d'éléments de $A$. Si pour tout $\mathcal{P}$ non ramifié dans $\text{Spec}(A)$, nous avons*

$$U_p \equiv V_p, \quad [\mathcal{P}],$$

*avec $p = N(\mathcal{P})$. Alors nous aurons $(U) = (V)$, c'est à dire $U_n = V_n$ pour tout $n \in \mathbb{N}$.*

**Démonstration.** La suite des différences : $(W) := (U) - (V)$ est encore de classe $C$ et vérifie $W_p \equiv 0, \quad [\mathcal{P}]$, pour tout $\mathcal{P}$ non ramifié. D'après le Théorème 18 et ses notations, nous aurons :

$$\begin{cases} 0 \equiv W_p \equiv W_1 & [\mathcal{P}], \ \forall \mathcal{P} \in \text{Split}(E/K) \\ 0 \equiv W_p \equiv s\, W_0 - W_1 & [\mathcal{P}], \ \forall \mathcal{P} \in \mathbb{P}(E/K, \tau) \end{cases}.$$

Puisque les ensembles $\text{Split}(E/K)$ et $\mathbb{P}(E/K, \tau)$ sont infinis, la première condition impose que nous ayons $W_1 = 0$ dans $A$. Et la seconde nous donne $s\, W_0 = 0$ dans $A$. Et donc que $W_0 = 0$ puisque $s \neq 0$. Comme les conditions initiales de $(W)$ sont nulles, cette suite est identiquement nulle. $\square$

### 4.3 Des congruences pour la fonction $\tau$ de Ramanujan

L'inspirante lecture de [7] que nous rencontrerons encore à la section suivante nous conduisit à trouver d'autres congruences de la fonction $\tau$ de Ramanujan dont nous ignorons l'originalité. Soit $l \geqslant 2$ un nombre premier. Pour tout $n \in \mathbb{N}$, la relation classique de Ramanujan-Mordell :

$$\tau(l^{n+2}) = \tau(l)\,\tau(l^{n+1}) - l^{11}\,\tau(l^n), \tag{7}$$

signifie que la suite $U: \mathbb{N} \to \mathbb{Z}, n \mapsto \tau(l^n)$, satisfasse la récurrence linéaire d'ordre deux :

$$U_{n+2} - \tau(l)\, U_{n+1} + l^{11}\, U_n = 0.$$

Le discriminant $\Delta(l)$ de son équation caractéristique vaut $\Delta(l) = \tau^2(l) - 4\, l^{11}$. Il est strictement négatif d'après le théorème de Deligne (anciennement conjecture de Ramanujan Petersson). En particulier $\Delta(l)$ n'est pas un carré dans $\mathbb{Z}$. On est donc dans le contexte d'application du théorème précédent. Par ailleurs, $U_0 = \tau(1) = 1$ et $U_1 = \tau(l) = s$ avec les notations précédentes. Ainsi $s\, U_0 - U_1 = 0$. On a donc démontré le résultat suivant :

**Corollaire 20.** *Soit $l \geqslant 2$ un nombre premier. Pour tout nombre premier $p$ premier à $\Delta(l) = \tau^2(l) - 4\, l^{11}$, on a les congruences suivantes satisfaites par la fonction $\tau$ de Ramanujan :*

- *Si $\left(\frac{\Delta(l)}{p}\right) = 1$ alors $\tau(l^p) \equiv \tau(l) \quad [p]$.*

- *Si $\left(\frac{\Delta(l)}{p}\right) = -1$ alors $\tau(l^p) \equiv 0 \quad [p]$.*

Afin de faire quelques vérifications numériques notamment en calculant des symboles de Legendre, je me suis aperçu que pour les petites valeurs de $l$, $\Delta(l)$ contenait beaucoup de petits facteurs carrés et que les nombres premiers sans facteurs carrés intervenant dans sa décomposition étaient tous comparativement plus grands que $l$. J'ignore si ce phénomène possède une explication accessible. En tout cas c'est ce que l'on constate dans la table suivante

| $l$ | 2 | 3 | 5 | 7 |
|---|---|---|---|---|
| $-\Delta(l)$ | $2^6 \times 119$ | $2^2\, 3^4 \times 11 \times 181$ | $2^4\, 5^2 \times 89 \times 4831$ | $2^2\, 3^2\, 7^2 \times 19 \times 29 \times 47 \times 167$ |

Cette table donne l'impression que $\Delta(l)$ est toujours divisible par $l^2$. C'est faux conformément au contre-exemple, $\Delta(23) \equiv 1[23]$. En complément du corollaire, on obtient des congruences fortes modulo $l$. En effet (7) donne l'implication suivante :

$$\tau(l^{n+1}) \equiv \tau(l)\,\tau(l^n) \quad [l^{11}] \Rightarrow \tau(l^n) \equiv \tau^n(l) \quad [l^{11}].$$

## 4.4 Le cas des extensions cyclotomiques

Soit $M \geqslant 2$ un entier et $E = \mathbb{Q}_M := \mathbb{Q}(\zeta)$, le corps des racines $M^{\text{ème}}$ de l'unité. On s'intéresse ici aux suites $(U)$ d'entiers, $(A = \mathbb{Z})$ de classes $C(X) = \Phi_M(X)$, où $\Phi_M$ est le $M^{\text{ème}}$ polynôme cyclotomique. Toutes ces suites sont d'ordre $d = \varphi(M) = \deg(\Phi_M(X))$. L'un de leurs attraits est d'être périodiques et de permettre des calculs élémentaires des expressions $U_p \mod [p]$ en fonction des conditions initiales $(U_0, ..., U_{d-1})$. Par ailleurs les suites de Lucas associées sont liées à la fonction $\mu$ de Mœbuis.

**Proposition 21.** *Les suites d'entiers de classe $\Phi_M(X)$ ont les propriétés suivantes :*

1. *Elles sont périodiques de période $M$.*

2. *Pour tout $n \in \mathbb{N}$, si $r$ est le reste de la division Euclidienne de $n$ par $M$, et si $R_r(X)$ est le reste de la division de $X^r$ par $\Phi_M(X)$, alors pour toute suite $(U)$ de classe $\Phi_M(X)$ nous aurons :*

$$U_n = U_r = R_r\{U_0\}.$$

**Démonstration.** Nous avons les congruences de polynômes :

$$X^M \equiv 1 \quad [\Phi_M(X)] \Rightarrow X^{M+n} \equiv X^n \quad [\Phi_M(X)]$$

donc d'après la Proposition 6,

$$U_M = U_0 \quad \text{et} \quad U_{M+n} = U_n.$$

ceci prouve le point 1. Pour le second, avec les notations de l'énoncé, nous avons

$$(X^n \equiv X^r \equiv R_r(X) \quad [\Phi_M(X)]) \quad \Rightarrow \quad (U_n = U_r = R_r\{U_0\}),$$

toujours d'après la Proposition 6. □

Par exemple pour $M = 12, d = \varphi(12) = 4, \Phi_{12}(X) = \Phi_6(X^2) = X^4 - X^2 + 1$. Si on s'intéresse aux $U_p \mod [p]$, avec $p$ premier, différent de 2 et 3, alors tous les restes $r$ seront inversibles mod $[12]$ et nous aurons :

| $R_1(X) = X$ | $p \equiv 1 \quad [12]$ | $U_p = U_1$ |
|---|---|---|
| $R_5(X) = X^3 - X$ | $p \equiv 5 \quad [12]$ | $U_p = U_3 - U_1$ |
| $R_7(X) = -X$ | $p \equiv 7 \quad [12]$ | $U_p = -U_1$ |
| $R_{11}(X) = X - X^3$ | $p \equiv 11 \quad [12]$ | $U_p = U_1 - U_3$ |

Pour l'étude des suites Lucas associées, notons $Z_M$ l'ensemble des racines $M^{\text{ème}}$ de l'unité et $\text{Prim}(M) \subset Z_M$, le sous ensemble des racines primitives. Puisque $\text{Prim}(M)$ coïncide avec l'ensemble des racines de $\Phi_M(X)$, la suite de Lucas de classe $\Phi_M(X)$ sera donnée par :

$$L_n^{(M)} := \sum_{z \in \text{Prim}(M)} z^n.$$

**Proposition 22.** *Les termes de la suite de Lucas se calculent conformément aux propriétés suivantes :*

1. *Pour $n=0$, on a $L_0^{(M)} = \varphi(M)$.*
2. *Si $n$ et $M$ sont premiers entre eux alors $L_n^{(M)} = L_1^{(M)}$.*
3. *Si $\mu$ désigne la fonction de Mœbuis alors $L_1^{(M)} = \mu(M)$.*
4. *En contraste du second point, si $n$ divise $M$ alors $L_n^{(M)} = \frac{\varphi(M)}{\varphi(M/n)} L_1^{(M/n)} = \frac{\varphi(M)}{\varphi(M/n)} \mu(M/n)$.*
5. *La règle du pgcd : Si $d := \text{pgcd}\{n; M\}$ alors : $L_n^{(M)} = \frac{\varphi(M)}{\varphi(M/d)} \mu(M/d) = L_d^{(M)}$.*

Ces règles sont utiles pour calculer les premiers termes de la suite de Lucas. Elle montrent que ces termes se réduisent à certaines valeurs de $L_1$; elles mêmes liées à la fonction de Mœbuis. Ensuite, la relation de récurrence permet de faire le reste. A titre d'illustration je donne les premières valeurs pour $M = 12$. Ici la récurrence est $L_{n+4}^{(12)} - L_{n+2}^{(12)} + L_n^{(12)} = 0 \quad (\mathcal{R})$.

| $n$ | 0 | 1 | 2 | 3 | 4 | 5 | 6 | 7 | 8 | 9 | 10 | 11 |
|---|---|---|---|---|---|---|---|---|---|---|---|---|
| $L_n^{(12)}$ | 4 | 0 | 2 | 0 | $-2$ | 0 | $-4$ | 0 | $-2$ | 0 | 2 | 0 |
| | R1 | R3 | R4 | $(\mathcal{R})$ | $(\mathcal{R})$ | R3,2 | $(\mathcal{R})$ | R3,2 | $(\mathcal{R})$ | $(\mathcal{R})$ | $(\mathcal{R})$ | R3,2 |

Dans la dernière ligne nous avons indiqué le numéro de la règle permettant de faire le calcul. Le symbole $(\mathcal{R})$ désigne la récurrence $(\mathcal{R})$. Par exemple $L_5 = L_3 - L_1 \Leftrightarrow 0 = L_3 - 0$, permet d'obtenir $L_3$ qui pourrait aussi se calculer grâce à R4. Cette règle R4 est avec R3 nécessaire au calcul de $L_2$ :

$$L_2^{(12)} = \frac{\varphi(12)}{\varphi(6)} L_1^{(6)} = \frac{4}{2} \mu(6) = 2 \mu(2 \times 3) = 2 \times (-1)^2 = 2.$$

**Démonstration.** Pour les Points 1 et 2. L'égalité $L_0^{(M)} = \varphi(M)$, provient du fait que $\text{card}(\text{Prim}(M)) = \varphi(M)$. La relation $L_n^{(M)} = L_1^{(M)}$, si $n$ et $M$ sont premiers entre eux provient du fait que dans ce contexte, $z \mapsto z^n$ induit une bijection comme application de $Z_M \to Z_M$ ou de $\text{Prim}(M) \to \text{Prim}(M)$.

Pour le Point 3 : Puisque $L_1^{(M)}$ est la somme des racines primitives $M^{\text{ème}}$ de l'unité, c'est au signe près le coefficient de $X^{\varphi(M)-1}$ du polynôme cyclotomique. C'est à dire :

$$\Phi_M(X) = X^{\varphi(M)} - L_1^{(M)} X^{\varphi(M)-1} + \ldots$$

Ainsi,

a) Quand $M = p$ est premier, $L_1^{(p)} = -1 = \mu(p)$ car $\Phi_p(X) = X^{p-1} + X^{p-2} + \cdots + X + 1$

b) Quand $M = p^s$ avec $s \geqslant 2$, est la puissance d'un nombre premier on a classiquement que $\Phi_{p^s}(X) = \Phi_p(X^{p^{s-1}})$. Donc $L_1^{(p^s)} = 0 = \mu(p^s)$.

c) Il reste à voir qu'à l'instar de la fonction de Mœbuis, $M \mapsto L_1^{(M)}$, est faiblement multiplicative. En effet, si $M$ et $M'$ sont premiers entre eux, l'application :

$$\text{Prim}(M) \times \text{Prim}(M') \to \text{Prim}(MM'), (z, z') \mapsto z \times z',$$

est bijective. Ainsi la multiplication des deux sommes $L_1^{(M)}$ et $L_1^{(M')}$ donnera $L_1^{(M \times M')}$.

Preuve du Point 4. Si $n$ divise $M$, considérons le diagramme commutatif :

$$\begin{array}{ccccccc}
(\mathbb{Z}/M\mathbb{Z})^* & \to & \operatorname{Prim}(M) & & \bar{k} & \mapsto & \exp\left(i\frac{2\pi\bar{k}}{M}\right)=z \\
\downarrow\Psi & & \downarrow z\mapsto z^n & , & \downarrow & & \downarrow \\
\left(\mathbb{Z}/\frac{M}{n}\mathbb{Z}\right)^* & \to & \operatorname{Prim}\left(\frac{M}{n}\right) & & \hat{k} & \mapsto & \exp\left(i\frac{2\pi\hat{k}}{M/n}\right)=z^n
\end{array}$$

Dans ce diagramme, les flèches horizontales sont des bijections. La flèches $\Psi$ n'est autre que la restriction du morphisme surjectif d'anneau : $\Psi\colon(\mathbb{Z}/M\mathbb{Z})\to\left(\mathbb{Z}/\frac{M}{n}\mathbb{Z}\right), \bar{k}\mapsto\hat{k}$. Ainsi quand on écrit la formule :

$$L_n^{(M)} = \sum_{z\in\operatorname{Prim}(M)} z^n = \sum_{\bar{k}\in(\mathbb{Z}/M\mathbb{Z})^*} \exp\left(i\frac{2\pi\bar{k}}{M/n}\right),$$

les $\varphi(M)$ exponentielles qui apparaissent dans la somme de droite sont toutes des racines primitives $(M/n)^{\text{ème}}$ de l'unité. Ce sont les mêmes racines qui apparaissent dans la somme des $\varphi(M/n)$ termes :

$$L_1^{(M/n)} = \sum_{\hat{k}\in\left(\mathbb{Z}/\frac{M}{n}\mathbb{Z}\right)^*} \exp\left(i\frac{2\pi\hat{k}}{M/n}\right),$$

Donc le rapport des deux sommes (quand elles ne sont pas simultanément nulles) n'est autre que le quotient de leurs nombres de termes autrement dit :

$$\frac{L_n^{(M)}}{L_1^{(M/n)}} = \frac{\varphi(M)}{\varphi(M/n)} \Leftrightarrow L_n^{(M)} = \frac{\varphi(M)}{\varphi(M/n)} L_1^{(M/n)} = \frac{\varphi(M)}{\varphi(M/n)} \mu(M/n).$$

Enfin, la règle 3 donne la dernière égalité.

Preuve du Point 5. Ecrivons $n=dn'$ et $M=dM'$ avec $\operatorname{pgcd}\{n';M'\}=1$. On a

$$L_n^{(M)} = \sum_{z\in\operatorname{Prim}(M)} z^{dn'} = \sum_{z\in\operatorname{Prim}(M)} (z^d)^{n'},$$

mais l'application $z\mapsto z^d$ de $\operatorname{Prim}(M)\to\operatorname{Prim}(M'=M/d)$ est une surjection d'un ensemble à $\varphi(M)$ éléments vers un ensemble à $\varphi(M/d)$ éléments. Ainsi nous avons

$$L_n^{(M)} = \frac{\varphi(M)}{\varphi(M/d)} L_{n'}^{(M')}.$$

Puisque $\operatorname{pgcd}\{n';M'\}=1$, $L_{n'}^{(M')} = L_1^{(M')} = \mu(M') = \mu(M/d)$, d'après les Points 2 et 3. Et l'on conclut d'après le Point 4. $\square$

# 5 Les cas non abeliens

## 5.1 Le résultat général

Reprenons les notations de la Section 1. Soit $C(X)$ unitaire dans $A[X]$ de racines simples $\{x_1,...,x_d\}$ dans $B$ la clôture intégrale de $A$ dans le corps de décomposition de $C(X)$. Soit $G(E/K)$ le groupe de Galois de $C(X)$. Pour tout $n\in\mathbb{N}$, soit $R_n(X)$ le reste de $X^n$ par $C(X)$. Ces polynômes sont à coefficients dans $A$ et admettent des expressions de la forme :

$$R_n(X) = \sum_{k=0}^{d-1} r_k(n) X^k = r_{d-1}(n) X^{d-1} + \cdots + r_1(n) X + r_0(n).$$

Pour toute suite $(U)$ de classe $C(X)$ dans $A$, d'après la Proposition 5, nous aurons :

$$U_n = R_n\{U_0\} = r_{d-1}(n) U_{d-1} + \cdots + r_1(n) U_1 + r_0(n) U_0.$$

**Les polynômes des classes de conjugaison :**

**Définition 23.** *Si pour tout $s \in G(E/K)$, $\Psi_s(X) = \sum_{k=0}^{d-1} \Psi_k(s) X^k$, désigne le polynôme permutant $\Psi_s(X)$ de la Définition 12, alors, pour chaque classe de conjuguais-on $\Lambda$ de $G(E/K)$, on associera les d polynômes suivants : Pour tout entier $0 \leq k \leq d-1$, nous poserons*

$$\Omega_{\Lambda;k}(R) = \prod_{s \in \Lambda} (R - \Psi_k(s)),$$

*Notons aussi : $\mathbb{P}(\Lambda)$, l'ensemble des $\mathcal{P}$ non ramifiés dans $\mathrm{Spec}(A)$ tels qu'il existe $Q$ divisant $\mathcal{P}$ dans $B$ avec $\mathrm{Frob}_Q \in \Lambda$.*

Avec ces notations nous obtenons le :

**Théorème 24.** *Soit $\Delta = \Delta_C$ le discriminant de $C(X)$. Pour toute classe de conjugaison $\Lambda$ de $G(E/K)$ et pour tout entier $0 \leq k \leq d-1$,*

1. *Le polynôme $\Omega_{\Lambda;k}(R)$ est de degré $\lambda = \mathrm{card}(\Lambda)$. Ses coefficients sont dans l'anneau $A[\frac{1}{\Delta}]$ Plus précisément $\Delta^\lambda \Omega_{\Lambda;k}(R) \in A[R]$.*

2. *Comme élément de $K[R]$, $\Omega_{\Lambda;k}(R)$ est soit irréductible soit une puissance d'un polynôme irréductible.*

3. *Pour tout $\mathcal{P} \in \mathbb{P}(\Lambda)$, si $p = N(\mathcal{P})$, les coefficients de l'égalité $U_p = \sum_{k=0}^{d-1} r_k(p) U_k$, vérifient chacun la congruence :*

$$\Delta^\lambda \Omega_{\Lambda;k}(r_k(p)) \equiv 0 \quad [\mathcal{P}].$$

Avant la preuve commençons par observer que les Points 3 et 4 du Théorème 4, sont des cas particuliers de ce résultat. En effet, si $\Lambda = \{\sigma\}$ est la classe de conjugaison d'un élément du centre de $G(E/K)$, alors chaque $\Omega_{\Lambda;k}(R) = R - \Psi_k(\sigma)$ est de degré $\lambda = 1$. Et les congruences du théorème s'écrivent

$$\Delta\, r_k(p) \equiv \Delta\, \Psi_k(\sigma) \quad [\mathcal{P}]$$

pour tout $\mathcal{P} \in \mathbb{P}(\{\sigma\})$. Comme chaque membre de ces congruences est dans $A$, en les multipliant par les $U_k$ et en sommant nous ré-obtenons :

$$\Delta \sum_{k=0}^{d-1} r_k(p) U_k \equiv \Delta \sum_{k=0}^{d-1} \Psi_k(\sigma) U_k \quad [\mathcal{P}]$$

$$\Delta\, U_p \equiv \Delta\, \Psi_\sigma\{U_0\} \quad [\mathcal{P}].$$

De la même façon, de $\Psi_{\mathrm{Id}}(X) = X = \sum_{k=0}^{d-1} \Psi_k(\mathrm{Id})\, X^k$, nous obtenons que pour $\Lambda = \{\mathrm{Id}\}$ et $\mathcal{P} \in \mathrm{Split}(E/K) = \mathbb{P}(\{\mathrm{Id}\})$, la congruence précédente devient

$$U_p \equiv U_1 \quad [\mathcal{P}].$$

**Démonstration.**
  **Point 1**. La relation $\Omega_{\Lambda;k}(R) = \prod_{s \in \Lambda} (R - \Psi_k(s))$, prouve que $\Omega_{\Lambda;k}(R)$ est de degré $\lambda = \mathrm{card}(\Lambda)$. D'après le Point 1 du Théorème 4, chaque racine $\Psi_k(s)$ de $\Omega_{\Lambda;k}(R)$ est dans le module $\frac{1}{\Delta} B$. Donc $\Delta^\lambda \Omega_{\Lambda;k}(R) \in B[R]$. Maintenant si nous faisons agir $G(E/K)$ sur les coefficients de ce polynôme, le Lemme 14, nous donne que pour tout $\theta \in G(E/K)$,

$$\begin{aligned}
\theta \cdot (\Omega_{\Lambda;k}(R)) &= \prod_{s \in \Lambda} (R - \theta(\Psi_k(s))) \\
&= \prod_{s \in \Lambda} (R - \Psi_k(\theta s \theta^{-1})) \\
\theta \cdot (\Omega_{\Lambda;k}(R)) &= \Omega_{\Lambda;k}(R).
\end{aligned}$$

Cette formule prouve que $\Omega_{\Lambda;k}(R) \in K[R]$ et plus précisément que

$$\Delta^\lambda\, \Omega_{\Lambda;k}(R) \in B[R] \cap K[R] = A[R].$$

**Point 2.** Par ailleurs, si $s \in \Lambda$, l'orbite des conjugués $\theta\, s\, \theta^{-1}$ pour $\theta$ décrivant $G(E/K)$ coïncide avec $\Lambda$. Donc $G(E/K)$ agit transitivement sur les racines de $\Omega_{\Lambda;k}(R)$. Ceci prouve que si ce polynôme n'est pas irréductible comme élément de $K[R]$, il coïncide avec une puissance d'un polynôme irréductible.

**Point 3.** Soit $(U)$ une suite de classe $C(X)$ dans $A$, soit $\mathcal{P} \in \mathbb{P}(\Lambda)$, $s \in \Lambda$ et $Q \in \mathrm{Spec}(B)$ divisant $\mathcal{P}$ tel que $s = \mathrm{Frob}_Q$. D'après (6), nous aurons les congruences

$$\Delta\, U_p = \Delta\, U_{N(\mathcal{P})} \equiv \Delta\, \Psi_s\{U_0\} \quad [Q]$$

$$\Delta\, U_p \equiv \Delta \sum_{k=0}^{d-1} \Psi_k(s)\, U_k \quad [Q]$$

$$\sum_{k=0}^{d-1} \Delta\, r_k(p)\, U_k \equiv \sum_{k=0}^{d-1} \Delta\Psi_k(s)\, U_k \quad [Q]$$

Comme cette dernière congruence est vraie pour toutes les conditions initiales $(U_0, ..., U_{d-1})$, ceci implique que nous ayons $\Delta\, r_k(p) \equiv \Delta\Psi_k(s) \quad [Q]$. En faisant agir un élément $\theta \in G(E/K)$ sur cette congruence et en tenant compte du fait que $\Delta$ et $r_k(p)$ sont invariants, nous obtenons :

$$\Delta\, r_k(p) \equiv \Delta\Psi_k(s) \quad [Q]$$
$$\Delta\, r_k(p) \equiv \Delta\, (\theta(\Psi_k(s))) \quad [\theta \cdot Q]$$
$$\Delta\, r_k(p) \equiv \Delta\Psi_k(\theta\, s\, \theta^{-1}) \quad [\theta \cdot Q]$$
$$\Delta(r_k(p) - \Psi_k(\theta\, s\, \theta^{-1})) \equiv 0 \quad [\theta \cdot Q]$$

En multipliant toutes ces relations pour des $\theta$ tels que $\theta\, s\, \theta^{-1}$ décrive $\Lambda$, nous obtenons que $\Delta^\lambda\, \Omega_{\Lambda;k}(r_k(p))$ est dans $Q$. Mais comme $\Delta^\lambda\, \Omega_{\Lambda;k}(r_k(p))$ appartient aussi à $A$, il est dans $\mathcal{P} = A \cap Q$. $\square$

## 5.2 Les récurrences d'ordre trois de groupe $\mathfrak{S}_3$

Le théorème précédent s'illustre bien dans le contexte des récurrences d'ordre trois, car les polynôme $\Omega(R)$ deviennent explicitement calculable et $\mathfrak{S}_3$ est le plus petit groupe non abelien. Soit $C(X) = X^3 + u\, X + v \in A[X]$ irréductible de groupe de Galois $\mathfrak{S}_3$. C'est à dire que son discriminant :

$$\Delta = -4\, u^3 - 27\, v^2, \tag{8}$$

n'est pas un carré dans $A$. Le groupe $\mathfrak{S}_3$, contient trois classes de conjugaisons $\Lambda_f$ pour chaque valeur $f \in \{1; 2; 3\}$, où l'entier $f$ désigne l'ordre des éléments de la classe. Précisément :

- $\Lambda_1 = \{\mathrm{Id}\}$.
- $\Lambda_3 = \{c := (1; 2; 3), c^2 = (3; 1; 2)\}$ est la classe des éléments d'ordre trois de $\mathfrak{S}_3$. Elle contient deux éléments. Ainsi : $\lambda_3 := \mathrm{card}(\Lambda_3) = 2$.
- $\Lambda_2 = \{\tau_3 := (1; 2), \tau_2 := (1; 3), \tau_1 := (2, 3)\}$ est la classe des transpositions de $\mathfrak{S}_3$. Elle contient trois éléments. Ainsi : $\lambda_2 := \mathrm{card}(\Lambda_2) = 3$.

En conformité avec ces notations nous poserons $\mathbb{P}_f := \mathbb{P}(\Lambda_f)$; c'est l'ensemble des $\mathcal{P}$ non ramifiés dans $A$ possédant un diviseur $Q$ dont le Frobenuis est dans $\Lambda_f$ c'est à dire est d'ordre $f$. Le théorème de Chebotarev prévoit alors que les densités respectives de ces ensembles soient :

$$d(\mathbb{P}_1) = \frac{1}{6}, \quad d(\mathbb{P}_2) = \frac{3}{6} = \frac{1}{2}, \quad d(\mathbb{P}_3) = \frac{2}{6} = \frac{1}{3}.$$

En revanche, les idéaux premiers qui les composent dépendent essentiellement de $C(X)$ et sont liées aux suites $(U)$ de classe $C$, c'est à dire vérifiant la récurrence linéaire d'ordre trois :

$$U_{n+3} + u\, U_{n+1} + v\, U_n = 0 \tag{9}$$

Simplifions l'écriture des restes en les écrivant :

$$R_n(X) = r_2(n)\, X^2 + r_1(n)\, X + r_0(n) = a_n X^2 + b_n X + c_n.$$

Alors la Proposition 5, donne $U_n = a_n U_2 + b_n U_1 + c_n U_0$ pour toute suite de classe $C(X)$ c'est à dire vérifiant (9). Avec ces notations nous obtenons le

**Théorème 25.** *Pour toute suite linéaire de classe $C(X)$ dans $A$, c'est à dire satisfaisant (9). Pour tout $\mathcal{P}$ non ramifié, si nous notons $p = N(\mathcal{P})$, alors on a une congruence linéaire*

$$U_p \equiv a\, U_2 + b\, U_1 + c\, U_0 \quad [\mathcal{P}],$$

*où les coefficients appartiennent à $A$, dépendent de l'idéal premier $\mathcal{P}$, mais vérifient des congruences polynomiales qui elles ne dépendent pas de $\mathcal{P}$.*

- *On aura toujours la « congruence de Lucas » : $3\,c \equiv 2\,u\,a \quad [\mathcal{P}]$.*
- *Si $\mathcal{P} \in \mathbb{P}_1$ alors $U_p \equiv U_1 \quad [\mathcal{P}]$*
- *Si $\mathcal{P} \in \mathbb{P}_3$ alors $\begin{cases} \Delta\,a^2 - 9\,u^2 & \equiv 0\ [\mathcal{P}] \\ \Delta\,(b^2 + b) - (u^3 + 27\,v^2) & \equiv 0\ [\mathcal{P}] \\ \Delta\,c^2 - 4\,u^4 & \equiv 0\ [\mathcal{P}] \end{cases}$.*
- *Si $\mathcal{P} \in \mathbb{P}_2$ alors $\begin{cases} \Delta\,a^3 - 9\,u^2\,a - 27\,v & \equiv 0\ [\mathcal{P}] \\ \Delta\,b^3 + 3\,u^3\,b - u^3 & \equiv 0\ [\mathcal{P}] \\ \Delta\,c^3 - 4\,u^4\,c - 8\,u^3\,v & \equiv 0\ [\mathcal{P}] \end{cases}$.*

Outre la congruence de Lucas qui est « universelle », les autres congruences polynomiales sur les coefficients sont discriminantes dans la mesure où elles dépendent de la classe de conjugaison du Frobenuis, mais pas de l'idéal premier. Observons aussi qu'elles ne sont pas intégrales, il y a à chaque fois un facteur $\Delta$ devant la plus haute puissance de chaque relation. C'est ce que prévoyait le Théorème 24, avec toutefois ici des puissances de $\Delta$, moins importantes que celles prévues. En revanches, les relations polynomiales données quand $\mathcal{P} \in \mathbb{P}_f$ si $f = 2$ ou 3, sont les relations $\Omega(R) \equiv 0 \quad [\mathcal{P}]$ et elles sont bien de degré $\lambda_f$. Enfin, la relation $U_p \equiv U_1 \quad [\mathcal{P}]$ si $\mathcal{P} \in \mathbb{P}_1$ est celle déjà prévue par le Théorème 4. Elle ne fut ajoutée que pour la complétude de l'énoncé.

**Remarques sur la caractère discriminant des congruences du Théorème 25.**

D'un point de vue pratique, si $R_p(X) = a\,X^2 + b\,X + c$ et que l'on souhaite tester la classe $\mathbb{P}_f$ de l'idéal premier $\mathcal{P}$, il est nécessaire de vérifier le fait que les trois systèmes de congruences ne sont pas simultanément compatibles. Autrement dit que nous avons la propriété suivante :

**Proposition 26.** *Soit $\mathcal{P}$ un idéal premier non ramifié de $A$ c'est à dire tel que $\Delta \notin \mathcal{P}$, alors pour tout $(a,b,c) \in A^3$ et tout couple $(f, f') \in \{1,2,3\}^2$ avec $f \neq f'$ les congruences des systèmes $\mathbb{P}_f$ et $\mathbb{P}_{f'}$ sont simultanément impossibles.*

**Démonstration.** Soit $(a,\,b,\,c) \in A^3$, un triplet quelconque, nous allons tester les conditions qu'imposent le fait qu'il satisfasse les congruences de $\mathbb{P}_f$ et de $\mathbb{P}_{f'}$.

Pour $\mathbb{P}_3$ et $\mathbb{P}_2$ : Pour les équations en $a$ et $b$, nous avons les divisions euclidiennes :

$$\begin{cases} \Delta\,a^3 - 9\,u^2\,a - 27\,v = a\,(\Delta\,a^2 - 9\,u^2) - 27\,v \\ \Delta\,b^3 + 3\,u^3\,b - u^3 = (\Delta\,b^2 + \Delta\,b - (u^3 + 27\,v^2))\,(b-1) - (2\,u^3 + 27\,v^2) \end{cases}$$

Donc si les deux congruences en $a$ sont simultanément vérifiées, cela impose que nous ayons $27\,v \in \mathcal{P}$. De même, si les deux congruences en $b$ sont simultanément vérifiées, nous aurons $2\,u^3 + 27\,v^2 \in \mathcal{P}$. Ceci implique $2\,u^3 \in \mathcal{P}$, puisque $27\,v \in \mathcal{P}$. Donc au final $\Delta = -(4\,u^3 + 27\,v^2) \in \mathcal{P}$ ce que nous avons supposé être faux.

Pour $\mathbb{P}_2$ et $\mathbb{P}_1$ : Que $(a,b,c)$ satisfasse les congruences de $\mathbb{P}_1$ signifie que $(a,b,c) \equiv (0,1,0) \quad [\mathcal{P}]$. Les équations en $a$ et $c$ donnent alors :

$$\begin{cases} \Delta\,a^3 - 9\,u^2\,a - 27\,v \equiv 0\ [\mathcal{P}] \Rightarrow 27\,v \equiv 0\ [\mathcal{P}] \\ \Delta\,b^3 + 3\,u^3\,b - u^3 \equiv 0\ [\mathcal{P}] \Rightarrow \Delta + 2\,u^3 = -(2\,u^3 + 27\,v^2) \equiv 0\ [\mathcal{P}] \end{cases}.$$

Et l'on retrouve les deux congruences précédentes : $27\,v \in \mathcal{P}$ et $2\,u^3 + 27\,v^2 \in \mathcal{P}$ qui conduisaient à $\Delta = -(4\,u^3 + 27\,v^2) \in \mathcal{P}$.

Pour $\mathbb{P}_3$ et $\mathbb{P}_1$ : Que $(a,b,c)$ satisfasse les congruences de $\mathbb{P}_1$ signifie que $(a,b,c) \equiv (0,1,0) \quad [\mathcal{P}]$. Les équations en $a, b$ et $c$ donnent alors :

$$\begin{cases} \Delta\,a^2 - 9\,u^2 & \equiv 0\ [\mathcal{P}] \Rightarrow 9\,u^2 & \equiv 0\ [\mathcal{P}] \\ \Delta\,(b^2 + b) - (u^3 + 27\,v^2) & \equiv 0\ [\mathcal{P}] \Rightarrow -9\,u^3 - 3 \times 27\,v^2 & \equiv 0\ [\mathcal{P}] \\ \Delta\,c^2 - 4\,u^4 & \equiv 0\ [\mathcal{P}] \Rightarrow 4\,u^4 & \equiv 0\ [\mathcal{P}] \end{cases}.$$

Les deux premières congruences impliquent $3 \times 27\, v^2 = (9\, v)^2 \in \mathcal{P} \Rightarrow 9\, v \in \mathcal{P}$ car $\mathcal{P}$ est premier. Toujours pour la même raison, la troisième implique que nous ayons $2\, u^2 \in \mathcal{P}$. D'où finalement $\Delta = -(4\, u^3 + 27\, v^2) \in \mathcal{P}$ ce que nous avons supposé être faux. □

Nous ne savons pas encore prouver de généralisation de cette propriété dans le contexte général du Théorème 24.

### 5.3 Preuve du Théorème 25

Malgré la simplicité du groupe $\mathfrak{S}_3$ les calculs des polynômes $\Omega_{\Lambda;k}(R)$, sont assez volumineux mais peuvent être rendus abordables au moyen de quelques simplifications que nous présentons maintenant.

#### 5.3.1 Relations de Newton et de Lucas.

Ici $C(X) = X^3 + u\, X + v = (X - x_1)\,(X - x_2)\,(X - x_3)$. Pour tout élément $\varphi \in E(X)$, nous noterons :
$$\oplus_c \varphi := \varphi + c \cdot \varphi + c^2 \cdot \varphi,$$
la « somme circulaire » associée à $\varphi$, où $c$ désigne le cycle $c = (1; 2; 3) \in \mathfrak{S}_3$ et $c \cdot \varphi$ l'action du cycle sur $\varphi$. Par exemple avec cette notation les relations de Newton et de Lucas s'écrivent ici :

$$\begin{cases} 3 &= x_1^0 + x_2^0 + x_3^0 &= \oplus_c 1 &= L_0 \\ 0 &= x_1 + x_2 + x_3 &= \oplus_c x_1 &= L_1 \\ u &= x_1 x_2 + x_2 x_3 + x_3 x_1 &= \oplus_c x_1 x_2 & \\ -v &= x_1 x_2 x_3 & & \\ -2\,u &= x_1^2 + x_2^2 + x_3^2 &= \oplus_c x_1^2 &= L_2 \\ u^2 &= x_1^2 x_2^2 + x_2^2 x_3^2 + x_3^2 x_1^2 &= \oplus_c x_1^2 x_2^2 & \end{cases} \qquad (10)$$

Les deuxième et troisième relations donnent la suivante :
$$u = -x_1^2 + x_2\, x_3 = -x_2^2 + x_3\, x_1 = -x_3^2 + x_1\, x_2 \qquad (11)$$

A partir de là, la relation donnant $L_2$ s'obtient en sommant ces trois dernières relations. Puisque la suite de Lucas vérifie $L_{n+3} + u\, L_{n+1} + v\, L_n = 0$ c'est à dire la récurrence (9), on obtient le tableau de ses premières valeurs :

| $n$ | 0 | 1 | 2 | 3 | 4 | 5 | 6 | 7 |
|---|---|---|---|---|---|---|---|---|
| $L_n$ | 3 | 0 | $-2\,u$ | $-3\,v$ | $2\,u^2$ | $5\,u\,v$ | $-2\,u^3 + 3\,v^2$ | $-7\,u^2\,v$ |

ainsi on vérifie bien que la congruence de la Proposition 7 : $L_p \equiv L_1 = 0 \quad [p]$, fonctionne bien pour les petits indices premiers.

Plus intéressant pour notre propos est la relation universelle sur les coefficients des congruences de la Proposition 8 : Si pour toute suite de classe $C(X)$ dans $A$, nous avons
$$U_p = U_{N(\mathcal{P})} = a\, U_2 + b\, U_1 + c\, U_0,$$
On aura en particulier pour la suite de Lucas :
$$L_p = L_{N(\mathcal{P})} = a\, L_2 + b\, L_1 + c\, L_0 \equiv L_1 \quad [\mathcal{P}].$$
Ce qui compte tenu du tableau des conditions initiales s'écrit encore
$$-2\,u\,a + 3\,c \equiv 0 \quad [\mathcal{P}] \Leftrightarrow 3\,c \equiv 2\,u\,a \quad [\mathcal{P}] \qquad (12)$$

Cette relation est le premier point du Théorème 25. Elle nous permettra de déduire une relation sur $a$ d'une relation sur $c$ et inversement.

#### 5.3.2 Différentes expressions de $C'(x_1)$ et discriminant.

En dérivant puis évaluant en $X = x_1$ la relation $C(X) = X^3 + u\, X + v = (X - x_1)\,(X - x_2)\,(X - x_3)$, nous obtenons
$$C'(x_1) = (x_1 - x_2)\,(x_1 - x_3) = 3\, x_1^2 + u,$$

et ses deux autres permutations circulaires. De sorte que si nous posons :

$$\delta := (x_1 - x_2)(x_1 - x_3)(x_2 - x_3),$$

cette grandeur est une racine carrée de $\Delta = \Delta_C$, et on obtient les trois égalités circulaires :

$$\delta = C'(x_1)(x_2 - x_3) = C'(x_2)(x_3 - x_1) = C'(x_3)(x_1 - x_2), \tag{13}$$

car $\delta$ est invariant par l'action du cycle $c = (1;2;3)$. Observons que si nous multiplions ces trois relations, en tenant compte de (10), nous obtenons très rapidement l'expression classique de $\Delta$ :

$$\begin{aligned}
\delta^3 &= C'(x_1)C'(x_2)C'(x_3)(-\delta) \\
-\Delta &= C'(x_1)C'(x_2)C'(x_3) \\
&= (u + 3x_1^2)(u + 3x_2^2)(u + 3x_3^2) \\
&= u^3 + 3u^2 \times \oplus_c x_1^2 + 9u \times \oplus_c x_1^2 x_2^2 + 27(x_1 x_2 x_3)^2 \\
&= u^3 + 3u^2 \times (-2u) + 9u \times (u^2) + 27v^2 \\
-\Delta &= 4u^3 + 27v^2.
\end{aligned}$$

Observons que la deuxième ligne de ce calcul n'est autre qu'un cas particulier de la formule classique du discriminant : $\Delta_f = (-1)^{n(n-1)/2} \Pi_i f'(x_i)$ voir ([4] exercice 12 p.193).

Si $\delta$ est invariant par l'action des trois cycles, il ne l'est pas sous celle des transpositions. D'autres quantités semblables vont intervenir dans les calculs des $\Omega(R)$. Par exemple si on développe $\delta$ on obtient 8 termes dont deux se simplifient pour donner :

$$\begin{aligned}
\delta &= (x_1^2 x_2 + x_2^2 x_3 + x_3^2 x_1) - (x_1 x_2^2 + x_2 x_3^2 + x_3 x_1^2) \\
&= S' - S \\
S' &:= \oplus_c x_1^2 x_2 \\
S &:= \oplus_c x_1 x_2^2
\end{aligned}$$

Pour calculer $S$ et $S'$ qui sont homogènes de degré trois en les $x_i$, on multiplie les relations : $0 \times u = (\oplus_c x_1) \times (\oplus_c x_1 x_2)$ de (10) pour obtenir :

$$\begin{cases} S' - S = \delta \\ S' + S = 3v \end{cases} \Leftrightarrow S' = \frac{3v + \delta}{2}, \quad S = \frac{3v - \delta}{2} \tag{14}$$

**Notation 27.** *Dans cette section comme dans la suivante, nous adapterons les lettres des variables à des notations simplifiées pour les rendre plus suggestives. Ainsi, conformément à la substitution :*

$$r_2 U_2 + r_1 U_1 + r_0 U_0 = a U_2 + b U_1 + c U_0,$$

*nous remplacerons la variable « R », par les lettres grecques « $\alpha, \beta, \gamma$ », dans les polynômes $\Omega(R)$. Autrement dit pour toute classe de conjugaison $\Lambda$, nous écrirons :*

$$\Omega_{\Lambda;2}(R) = \Omega_{\Lambda;2}(\alpha), \quad \Omega_{\Lambda;1}(R) = \Omega_{\Lambda;1}(\beta), \quad \Omega_{\Lambda;0}(R) = \Omega_{\Lambda;0}(\gamma).$$

### 5.3.3 Les polynômes $\Omega_{\Lambda;k}(R)$ pour $\Lambda = \Lambda_3$ : la classe des 3-cycles

**Calcul de $\Psi_c(X)$**

D'après le Lemme 13, nous avons :

$$\begin{aligned}
\Psi_c(X) &= \bigoplus_c \frac{C(X) x_2}{(X - x_1) C'(x_1)} \\
&= \bigoplus_c \frac{(X - x_2)(X - x_3) x_2}{C'(x_1)} \\
\Psi_c(X) &= \bigoplus_c \frac{(X - x_2)(X - x_3) x_2 (x_2 - x_3)}{C'(x_1)(x_2 - x_3)} \\
\delta \Psi_c(X) &= \bigoplus_c x_2 (x_2 - x_3)(X - x_2)(X - x_3). \tag{15}
\end{aligned}$$

La dernière égalité étant une conséquence de (13). Il reste à développer ce polynôme de degré deux pour l'identifier à :

$$\delta \Psi_c(X) = \delta \Psi_0(c) + \delta \Psi_1(c) X + \delta \Psi_2(c) X^2$$

**Preuve que : $\Omega_{\Lambda_3;2}(\alpha) = \alpha^2 - \frac{9 u^2}{\Delta}$**

Le coefficient de degré deux de $\delta \Psi_c(X)$ vaut

$$\begin{aligned} \delta \Psi_2(c) &= \bigoplus_c x_2 (x_2 - x_3) \\ &= \bigoplus_c x_2^2 - \bigoplus_c x_2 x_3 \\ \delta \Psi_2(c) &= -2 u - u \\ \Psi_2(c) &= \frac{-3 u}{\delta} \end{aligned}$$

d'après (10). En faisant agir l'une des transpositions sur cette relation, le Lemme 14, donne

$$\begin{cases} \tau(\Psi_2(c)) = \Psi_2(\tau c \tau^{-1}) = \Psi_2(c^2) \\ \tau\left(\frac{-3 u}{\delta}\right) = \frac{3 u}{\delta} \end{cases} \Rightarrow \Psi_2(c^2) = \frac{3 u}{\delta}$$

Ainsi en partant de la définition :

$$\Omega_{\Lambda_3;2}(\alpha) = (\alpha - \Psi_2(c))(\alpha - \Psi_2(c^2)) = \left(\alpha + \frac{3 u}{\delta}\right)\left(\alpha - \frac{3 u}{\delta}\right) = \alpha^2 - \frac{9 u^2}{\Delta}.$$

**Preuve que : $\Omega_{\Lambda_3;0}(\gamma) = \gamma^2 - \frac{4 u^4}{\Delta}$**

Un calcul équivalent peut être fait pour $\Omega_{\Lambda_3;0}(\gamma)$, mais il est plus rapide de le déduire de la relation de Lucas $3 c \equiv 2 u a \quad [\mathcal{P}]$ qui impose juste le changement de variable : $3 \gamma = 2 u \alpha$. Puis, à un facteur de proportion près d'avoir $\Omega_{\Lambda_3;0}(\gamma) = \Omega_{\Lambda_3;2}(3 \gamma / 2 u)$.

**Preuve que : $\Omega_{\Lambda_3;1}(\beta) = \beta^2 + \beta - \frac{(27 v^2 + u^3)}{\Delta}$**

Toujours d'après (15), le coefficient de degré un vaut

$$\begin{aligned} \delta \Psi_1(c) &= \bigoplus_c x_2 (x_2 - x_3) x_1 \\ &= \bigoplus_c x_2^2 x_1 - \bigoplus_c x_2 x_3 x_1 \\ \delta \Psi_1(c) &= S + 3 v \\ &= \frac{3 v - \delta}{2} + 3 v \\ \Psi_1(c) &= \frac{9 v - \delta}{2 \delta} \end{aligned}$$

d'après (10) et (14). En faisant agir l'une des transpositions sur cette relation, le Lemme 14, donne

$$\tau(\Psi_1(c)) = \Psi_1(\tau c \tau^{-1}) = \Psi_1(c^2) = \frac{-(9 v + \delta)}{2 \delta}.$$

Ainsi en partant de la définition puis en développant en tenant compte de (8), nous obtenons :

$$\Omega_{\Lambda_3;1}(\beta) = (\beta - \Psi_1(c))(\beta - \Psi_1(c^2)) = \left(\beta + \frac{\delta - 9 v}{2 \delta}\right)\left(\beta + \frac{\delta + 9 v}{2 \delta}\right) = \beta^2 + \beta - \frac{(27 v^2 + u^3)}{\Delta}.$$

### 5.3.4 Les polynômes $\Omega_{\Lambda;k}(R)$ pour $\Lambda = \Lambda_2$ : la classe des transpositions

Puisque $\Lambda_2 = \{\tau = \tau_3 := (1; 2), \tau_2 := (1; 3), \tau_1 := (2, 3)\}$ contient trois éléments, les polynômes, $\Omega_{\Lambda_2;k}(R)$ seront de degré trois.

**Calcul de $\Psi_\tau(X)$**

D'après le Lemme 13,

$$\Psi_\tau(X) = \frac{C(X)\,x_2}{(X-x_1)\,C'(x_1)} + \frac{C(X)\,x_1}{(X-x_2)\,C'(x_2)} + \frac{C(X)\,x_3}{(X-x_3)\,C'(x_3)}$$

$$X = \Psi_{\mathrm{Id}}(X) = \frac{C(X)\,x_1}{(X-x_1)\,C'(x_1)} + \frac{C(X)\,x_2}{(X-x_2)\,C'(x_2)} + \frac{C(X)\,x_3}{(X-x_3)\,C'(x_3)}$$

$$\Psi_\tau(X) - X = (x_1-x_2)\left[\frac{C(X)}{(X-x_2)\,C'(x_2)} - \frac{C(X)}{(X-x_1)\,C'(x_1)}\right]$$

$$\Psi_\tau(X) - X = (x_1-x_2)\left[\frac{C(X)\,(x_3-x_1)}{(X-x_2)\,\delta} - \frac{C(X)\,(x_2-x_3)}{(X-x_1)\,\delta}\right]$$

$$\delta\,[X - \Psi_\tau(X)] = (x_1-x_2)(X-x_3)\,[(x_1-x_3)(X-x_1) + (x_2-x_3)(X-x_2)]$$

$$= (x_1-x_2)(X-x_3)\,[-3\,x_3\,X + x_3 x_1 - x_1^2 + x_3 x_2 - x_2^2]$$

$$= (x_1-x_2)(X-x_3)\,[-3\,x_3\,X + (x_3 x_1 - x_2^2) + (x_3 x_2 - x_1^2)]$$

$$\delta\,[X - \Psi_{\tau_3}(X)] = (x_1-x_2)(X-x_3)\,[-3\,x_3\,X + 2\,u] \qquad (16)$$

Pour le passage de la ligne 3 à 4, nous avons utilisé (13). La relation (11) a permis la simplification de la dernière ligne. Il reste à développer ce polynôme de degré deux pour identifier :

$$\delta\,\Psi_{\tau_3}(X) = \delta\,\Psi_0(\tau_3) + \delta\,\Psi_1(\tau_3)\,X + \delta\,\Psi_2(\tau_3)\,X^2$$

**Preuve que :** $\Omega_{\Lambda_2;0}(\gamma) = \gamma^3 - \frac{4\,u^4}{\Delta}\,\gamma - \frac{8\,u^3 v}{\Delta}$

Le coefficient de degré zéro de $\delta\,\Psi_{\tau_3}(X)$ vérifie :

$$-\delta\,\Psi_0(\tau_3) = (x_1-x_2)\,(-x_3)\,2\,u \Rightarrow y_3 := \Psi_0(\tau_3) = \frac{2\,u(x_1-x_2)x_3}{\delta}.$$

En faisant agir le cycle $c$ sur cette relation nous obtenons :

$$\begin{cases} y_3 = \frac{2\,u(x_1-x_2)x_3}{\delta} \\ y_1 = \frac{2\,u(x_2-x_3)x_1}{\delta} \\ y_2 = \frac{2\,u(x_3-x_1)x_2}{\delta} \end{cases}$$

Et comme $\Omega_{\Lambda_2;0}(\gamma) = (\gamma - y_3)\,(\gamma - y_1)\,(\gamma - y_2)$, les coefficients extrêmes sont aisés à obtenir :

$$\bigoplus_c y_3 = 0, \quad y_3\,y_1 y_2 = \frac{8\,u^3}{\delta^3}(-\delta)(-v) = \frac{8\,u^3 v}{\Delta}.$$

Pour le coefficient de $\gamma$ :

$$\bigoplus_c y_3\,y_1 = \frac{4\,u^2}{\Delta}\bigoplus_c (x_1-x_2)(x_2-x_3)x_3\,x_1$$

$$= \frac{4\,u^2}{\Delta}\bigoplus_c (-C'(x_2))x_3\,x_1$$

$$= -\frac{4\,u^2}{\Delta}\bigoplus_c (3\,x_2^2 + u)x_3\,x_1$$

$$= -\frac{4\,u^2}{\Delta}\bigoplus_c [-3\,v\,x_2 + u\,x_3 x_1]$$

$$\bigoplus_c y_3\,y_1 = -\frac{4\,u^4}{\Delta}$$

Et l'expression développée de $\Omega_{\Lambda_2;0}(\gamma) = \gamma^3 - \frac{4\,u^4}{\Delta}\,\gamma - \frac{8\,u^3 v}{\Delta}$, s'en déduit.

**Preuve que :** $\Omega_{\Lambda_2;2}(\alpha) = \alpha^3 - \frac{9\,u^2}{\Delta}\alpha - \frac{27\,v}{\Delta}$

Comme précédemment on passera du calcul de $\Omega_{\Lambda_2;0}(\gamma)$ à celui de $\Omega_{\Lambda_2;2}(\alpha)$ au moyen de la relation de Lucas.

**Preuve que : $\Omega_{\Lambda_2;1}(\beta) = \beta^3 + \frac{3\,u^3}{\Delta}\beta - \frac{u^3}{\Delta}$**

Isolons maintenant le terme en $X$ de (16) :

$$\begin{aligned}
\delta\,[1-\Psi_1(\tau_3)] &= (x_1-x_2)\,[2\,u+3\,x_3^2] \\
(x_1-x_2)\,C'(x_3) - \delta\Psi_1(\tau_3) &= (x_1-x_2)\,[2\,u+3\,x_3^2] \\
\delta\Psi_1(\tau_3) &= (x_1-x_2)\,[C'(x_3)-2\,u-3\,x_3^2] \\
z_3 := \Psi_1(\tau_3) &= \frac{-u\,(x_1-x_2)}{\delta}
\end{aligned}$$

Pour obtenir ces différentes expressions nous avons encore une fois utilisé les relations liant $C'(x_3)$ et $\delta$. Quand aux permutations circulaires de cette dernière relation elles nous donnent :

$$\begin{cases} z_3 = \frac{-u\,(x_1-x_2)}{\delta} \\ z_1 = \frac{-u\,(x_2-x_3)}{\delta} \\ z_2 = \frac{-u\,(x_3-x_1)}{\delta} \end{cases},$$

puis $\Omega_{\Lambda_2;1}(\beta) = (\beta-z_3)(\beta-z_1)(\beta-z_2)$. Comme précédemment les termes extrêmes se calculent vite :

$$\bigoplus_c z_3 = 0, \quad z_3\,z_1\,z_2 = \frac{-u^3}{\delta^3}(-\delta) = \frac{u^3}{\Delta}.$$

Par ailleurs,

$$z_3\,z_1 = \frac{-u^2}{\Delta}(x_2-x_1)(x_2-x_3) = \frac{-u^2}{\Delta}C'(x_2) = \frac{-u^2}{\Delta}[3\,x_2^2+u].$$

Ainsi,

$$\bigoplus_c z_3\,z_1 = \frac{-u^2}{\Delta}[3\times(-2\,u)+3\,u] = \frac{3\,u^3}{\Delta}.$$

Et le résultat final s'en déduit.

## 5.4 L'exemple de Wilton et Serre et les suites de Padovan et Lucas

Ce cas particulier nous fut inspiré par la lecture de [7]. Les suites linéaires d'ordre trois d'entiers relatifs qui ressemblent le plus à Fibonacci sont par exemple celles qui satisfont la récurrence :

$$U_{n+3} = U_{n+1} + U_n. \tag{17}$$

Dans [8] et son article : « La sculpture et les nombres », Ian Stewart évoque deux suites célèbres satisfaisant cette récurrence : la suite $(P)$ dite de Padovan et celle de Lucas $(\mathcal{L})$ de premiers termes respectifs :

| $n$ | 0 | 1 | 2 | 3 | 4 | 5 | 6 | 7 |
|---|---|---|---|---|---|---|---|---|
| $P_n$ | 1 | 1 | 1 | 2 | 2 | 3 | 4 | 5 |
| $\mathcal{L}_n$ | 3 | 0 | 2 | 3 | 2 | 5 | 5 | 7 |

La suite de Padovan permet le tracé d'une jolie spirale constituée de triangles équilatéraux, quant à Edouard Lucas il avait déjà observé en 1876 l'implication :

$$n \text{ premier} \Rightarrow \mathcal{L}_n \equiv 0 \quad [n].$$

Suite à ses travaux on a longtemps cherché des nombres « pseudo premiers » c'est à dire vérifiant $\mathcal{L}_n \equiv 0 \ [n]$ sans que $n$ soit premier. Ici contrairement à la suite de Lucas associée à la récurrence de Fibonacci il faut aller chercher les nombres « pseudo premiers » très loin. Le plus petit d'entre eux est $n = 541^2 = 271441$.

Le polynôme caractéristique des récurrences (17) est $C(X) = X^3 - X - 1 = X^3 + u\,X + v$. Il est irréductible, étant un polynôme d'Artin-Schreier. Son discriminant vaut $\Delta = -23$, ainsi son groupe de Galois sur $\mathbb{Q}$ est $\mathfrak{S}_3$. Ici, la spécialisation du Théorème 25, c'est à dire en posant $u = v = -1$, donne le :

**Corollaire 28.** *Pour toute suite linéaire de classe $C(X)$ dans $\mathbb{Z}$, c'est à dire satisfaisant $U_{n+3} = U_{n+1} + U_n$. Pour tout $p \neq 23$ non ramifié, on a une congruence linéaire :*

$$U_p = a\, U_2 + b\, U_1 + c\, U_0 \quad [p],$$

*où les coefficients appartiennent à $\mathbb{Z}$, dépendent de $p$, mais vérifient des congruences polynomiales qui elles ne dépendent pas de $p$.*

- *On aura toujours la « congruence de Lucas » : $3\,c \equiv -2\,a \quad [p]$.*
- *Si $p \in \mathbb{P}_1$ alors $U_p \equiv U_1 \quad [p]$*
- *Si $p \in \mathbb{P}_3$ alors $\begin{cases} 23\,a^2 + 9 & \equiv\ 0 \ [p] \\ 23\,(b^2 + b) + 26 & \equiv\ 0 \ [p] \\ 23\,c^2 + 4 & \equiv\ 0 \ [p] \end{cases}$.*
- *Si $p \in \mathbb{P}_2$ alors $\begin{cases} 23\,a^3 + 9\,a - 27 & \equiv\ 0 \ [p] \\ 23\,b^3 + 3\,b - 1 & \equiv\ 0 \ [p] \\ 23\,c^3 + 4\,c + 8 & \equiv\ 0 \ [p] \end{cases}$.*

Dans ce cas particulier, le corps de décomposition $E/\mathbb{Q}$ de $C(X)$ coïncide avec $K_{23}$ : la plus grande sous extension de $\bar{\mathbb{Q}}$ qui est non ramifié en dehors de 23. C'est le fait que $E/\mathbb{Q}(\sqrt{-23})$ soit le corps de classe absolu de $\mathbb{Q}(\sqrt{-23})$ qui permette d'obtenir dans ce contexte la description complète des ensembles $\mathbb{P}_f$ cité dans [7] :

**Théorème 29.** *Si $E/\mathbb{Q}$ est le corps de décomposition de $X^3 - X - 1$, alors pour tout nombre premier $p \neq 23$, on a :*

- *$p \in \mathbb{P}_1$ ssi $p$ est représenté par la forme quadratique $x^2 + 23\,y^2$*
- *$p \in \mathbb{P}_3$ ssi $\left(\frac{p}{23}\right) = 1$, mais $p$ n'est représenté par la forme quadratique $x^2 + 23\,y^2$.*
- *$p \in \mathbb{P}_2$ ssi $\left(\frac{p}{23}\right) = -1$.*

Le Corollaire 28, permet de trouver numériquement pour chaque nombre premier donné l'ensemble $\mathbb{P}_f$ auquel il appartient. En effet, pour chaque $p \neq 23$ on peut calculer sont reste $R_p(X) = a\,X^2 + b\,X + c$, et on regarde le groupe de congruences donné par le corollaire que satisfont les trois nombres $(a, b, c)$. C'est très pratique et aisé, car les groupes de congruences sont discriminants. Mais peut-on déduire le Théorème 29, directement du corollaire 28 ? Voilà qui semble très loin d'être possible ! Le Théorème 29 appartient à la théorie du corps de classe que nous n'abordons pas dans cette étude.

Dans une vidéo : Wiles's marvellous proof, Henri Diamon attribue une variante du Théorème 29 à un travail de Hecke sur les formes modulaires ce qui est aussi l'esprit de l'article de Serre déjà cité.

Signalons que la partie facile du Théorème 29, est contenue dans la :

**Proposition 30.** *Soit $C(X) = X^3 + u\,X + v \in A[X]$ de groupe $\mathfrak{S}_3$. Soit un idéal premier de $A$ non ramifié dans $B$. Alors $\mathcal{P} \in \mathbb{P}_2$ si et seulement si $\Delta_C$ n'est pas un carré dans $\mathbb{F}_\mathcal{P}$.*

Ainsi, même dans le cas où $A = \mathbb{Z}$, toute la difficulté est de savoir comment les nombres premiers tels que $\left(\frac{\Delta_C}{p}\right) = 1$, se répartissent entre $\mathbb{P}_1$ et $\mathbb{P}_3$. Ici nous gardons les notations de la Section 5.2.

**Démonstration.** Soit $\mathcal{P}$ non ramifié et $Q$ au dessus, soit $\sigma = \mathrm{Frob}_Q$ vu comme élément de $\mathrm{Gal}(\mathbb{F}_Q / \mathbb{F}_\mathcal{P})$ avec $\mathbb{F}_Q = \mathbb{F}_\mathcal{P}[\bar{x}_1, ..., \bar{x}_d]$ d'après le Lemme 16. Examinons les trois possibilités :

- Si $\mathcal{P} \in \mathbb{P}_1$, les $\bar{x}_i$ appartiennent à $\mathbb{F}_\mathcal{P}$ donc $\bar{\delta}$ aussi. Ainsi $\overline{\Delta_C} = \bar{\delta}^2$ est un carré dans $\mathbb{F}_\mathcal{P}$.
- Si $\mathcal{P} \in \mathbb{P}_3$, $\sigma$ agit comme un trois cycle, il laisse donc $\bar{\delta}$ invariant. Ainsi $\overline{\Delta_C} = \bar{\delta}^2$ est encore un carré dans $\mathbb{F}_\mathcal{P}$.
- Si $\mathcal{P} \in \mathbb{P}_2$, quitte à renuméroter les racines, on peut supposer que $\sigma$ fixe $\bar{x}_1$ et transpose $\bar{x}_2$ et $\bar{x}_3$. Cela implique que $\bar{x}_2$ et $\bar{x}_3$ soient quadratiques sur $\mathbb{F}_\mathcal{P}$ de minimal $f(X) := (X - \bar{x}_2)(X - \bar{x}_3)$. Or
$$\overline{\Delta_C} = f^2(\bar{x}_1) \times (\bar{x}_2 - \bar{x}_3)^2,$$
avec $f(\bar{x}_1) \in \mathbb{F}_\mathcal{P}$. Mais $(\bar{x}_2 - \bar{x}_3)^2$ coïncide avec le discriminant de $f(X)$ lequel n'est pas un carré car $(\bar{x}_2 - \bar{x}_3) \notin \mathbb{F}_\mathcal{P}$. Ainsi $\overline{\Delta_C}$ n'est pas un carré dans $\mathbb{F}_\mathcal{P}$ dans ce contexte. $\square$

### 5.5 La suite de « Cartier Trinks » et le groupe $\mathbb{PSL}(2, \mathbb{F}_7)$

On doit à Pierre Cartier dans les années 1970, le premier exemple de polynôme de groupe de Galois $\mathbb{PSL}(2, \mathbb{F}_7)$ sur $\mathbb{Q}$. Il s'agit du polynôme $T(X) := X^7 - 7X + 3$. Il n'est pas aisé de montrer que ce polynôme admet $\mathbb{PSL}(2, \mathbb{F}_7)$ comme groupe de Galois. C'est vraisemblablement la première fois que des mathématiciens durent avoir recours à des ordinateurs pour y parvenir. Dans cette veine, on peut déterminer $\mathbb{P}(\{\text{Id}\}) = \text{Split}(T(X))$ au moyen des récurrences linéaires d'équation caractéristique $T(X)$. Il n'y a pas beaucoup de nombre premier totalement décomposés puisqu'ici leur densité est $1/168$. Néanmoins les premières valeurs sont celles de la liste suivante :

| $n$   | 289  | 709  | 759  | 1098 | 1108 |
|-------|------|------|------|------|------|
| $p_n$ | 1879 | 5381 | 5783 | 8819 | 8893 |

Pour remplir ce tableau nous avons utilisé le fait que la suite des restes vérifiait la récurrence :

$$R_{n+7}(X) - 7R_{n+1}(X) + 3R_n(X) = 0.$$

Puis utilisé le critère : $R_p(X) \equiv X \quad [p]$, donné par le Théorème 15. Ces expériences sur ordinateur furent réalisées par Timothée Schmoderer que je remercie ici.

## 6 Algorithmes de Berlekamp et Cantor-Zassenhaus

Ces deux algorithmes [1] et [2] permettent de factoriser un polynôme dans les corps finis $\mathbb{F}_q[x]$, c'est à dire en particulier de trouver leurs racines dans $\mathbb{F}_q$. l'algorithme de Berlekamp est moins efficace que celui de Cantor-Zassenhaus, mais ce dernier ne fonctionne que si on sait à priori que tous les facteurs irréductibles seront de même degré. Il y a un lien entre le présent travail et ces deux algorithmes. Les calculs de restes $R_n(X)$ ne permettent pas de factoriser le polynôme $C(X)$ sur des corps finis mais de déterminer le type de la factorisation et plus précisément le type galoisien de l'automorphisme de Frobenuis, dans le groupe de Galois de $C(X)$. Dans l'exemple de Cartier-Trinks de groupe $G = \mathbb{PSL}_2(\mathbb{F}_7) \simeq \mathbb{PGL}_3(\mathbb{F}_2)$. L'isomorphisme $G \simeq \mathbb{PGL}_3(\mathbb{F}_2)$, permet de voir $G$ agir sur les 7 racines de $C(X)$ que l'on identifie avec les 7 points du plan de Fano $\mathbb{P}_2(\mathbb{F}_2)$ :

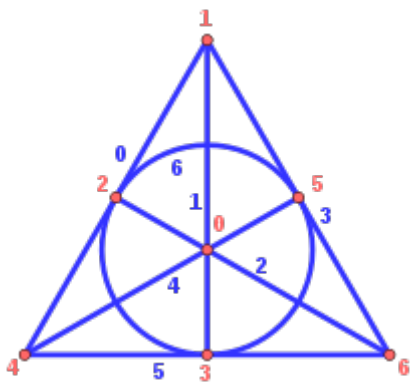

**Figure 1.**

Cette identification permet de classer les éléments de $G$ selon les types de permutations comme éléments de $\mathfrak{A}_7$ conformément au tableau suivant :

| Cas | type $\sigma$ | ordre de $\sigma$ | nombre de tels $\sigma$ |
|-----|---------------|-------------------|-------------------------|
| 1   | $1^7$         | 1                 | 1                       |
| 2   | $7 = 7^1$     | 7                 | 48                      |
| 3   | $2^2 \times 1^3$ | 2              | 21                      |
| 4   | $2^1 \times 4^1 \times 1^1$ | 4    | 42                      |
| 5   | $3^2 \times 1^1$ | 3              | 56                      |

Par exemple, dans le Cas 3, dire que $\sigma$ est de type $2^2 \times 1^3$ signifie qu'il agit sur les racines comme un produit de deux transpositions $(a,b)(c,d)$ et qu'il admet donc trois points fixes. Si un tel $\sigma$ coincide avec un $\text{Frob}_p$ cela signifiera que $C(X)$ se factorisera dans $\mathbb{F}_p[x]$ en produit de deux polynômes de degré deux et trois facteurs linéaires.

Dans le Cas 4, $\sigma$ est du type $\sigma = (a,b)\,(c,d,e,f)$ et si $\sigma = \text{Frob}_p$ alors dans $\mathbb{F}_p[x]$, $C(X)$ se factorisera en un terme de degré 2, un terme de degré 4 et un terme linéaire. Ainsi, ce que l'on voit ici c'est que pour ce qu'il en est de l'algorithme de Cantor-Zassenhaus, il ne peut fonctionner que dans les cas 1 et 2. Ce qui conformément au théorème de Chebotarev ne le rend opérant que pour

$$\frac{1}{168} + \frac{48}{168} = \frac{49}{168} \simeq 29/100,$$

des nombres premiers.

Notre travail permet théoriquement de déterminer le type du Frobenuis c'est à dire le type de la factorisation dans $\mathbb{F}_p[x]$ sans calculer cette dernière, mais juste en calculant le reste $R_p(X)$. Par exemple les restes $R_p(X)$ permettent de trouver que $p=1879$ est le plus petit premier où $C(X)$ se décompose totalement. Les algorithmes de Berlekamp et de Cantor-Zassenhaus permettent alors de trouver cette factorisation explicite dans $\mathbb{F}_{1879}[x]$ :

$$C(x) = (x+82)\,(x+298)\,(x+407)\,(x+883)\,(x+911)\,(x+1371)\,(x+1685).$$

Je remercie très chaleureusement Joris Van der Hoeven de m'avoir fait découvrir ces deux algorithmes. Dans [3], Van der Hoeven et ses collaborateurs étudient les complexités de ces algorithmes de factorisations sur les corps finis. L'étude de la complexité de notre travail reste à faire.

# Bibliographie